\documentclass[12pt,draft]{article}

\usepackage[cp1251]{inputenc}
\usepackage[T2A]{fontenc}
\usepackage[english,ukrainian]{babel}
\usepackage{amsmath,amsfonts,amssymb}
\usepackage{geometry}
\begin{document}

\begin{center}

\textbf{Anna Anop and Aleksandr Murach}

{\small(Institute of Mathematics of the National Academy of Sciences of Ukraine, Kyiv)}

\bigskip

\textbf{INTERPOLATION SPACES OF GENERALIZED SMOOTHNESS\\ AND THEIR APPLICATIONS TO ELLIPTIC EQUATIONS}

\bigskip\bigskip

\textbf{Анна Аноп і Олександр Мурач}

{\small(Інститут математики Національної академії наук України, Київ)}

\bigskip

\textbf{ІНТЕРПОЛЯЦІЙНІ ПРОСТОРИ УЗАГАЛЬНЕНОЇ ГЛАДКОСТІ\\ТА ЇХ ЗАСТОСУВАННЯ ДО ЕЛІПТИЧНИХ РІВНЯНЬ}

\bigskip

{\small 12 жовтня 2021 р., доповнено 1 червня 2023 р.}

\end{center}


\bigskip

\textbf{Abstract.} We introduce and investigate classes of normed or quasinormed distribution spaces of generalized smoothness that can be obtained by various interpolation methods applied to classical Sobolev, Nikolskii\,--\,Besov, and Triebel\,--\,Lizorkin spaces. An arbitrary positive function O-regularly varying at infinity serves as the order of regularity for the spaces introduced. They are broad generalizations of the above classical spaces and allow being well defined on smooth manifolds. We give applications of the spaces under investigation to elliptic equations and elliptic problems on smooth manifolds.

\bigskip

\textbf{1. Вступ.} Класичні шкали функціональних просторів Соболєва, Нікольського\,--\,Бєсова і Лізоркіна\,--\,Трібеля, градуйовані за допомогою числових параметрів, мають фундаментальне значення у сучасній математиці, перш за усе у математичному аналізі, теорії диференціальних рівнянь з частинними похідними, обчислювальній математиці (див., наприклад, монографії [1--10]). Вивчення та застосування цих просторів суттєво полегшується завдяки їх різним інтерполяційним властивостям. Так, простори Соболєва дробового порядку та простори Нікольського\,--\,Бєсова довільного порядку отримуються відповідно методами комплексної та дійсної інтерполяції просторів Соболєва цілого порядку. Крім того, вказані шкали просторів замкнені відносно різних методів інтерполяції  \cite[п.~2.4]{Triebel80}. Оскільки при інтерполяції успадковується обмеженість лінійних операторів, то багато властивостей цих просторів можна вивести з відповідних властивостей просторів цілого порядку, зокрема, просторів Соболєва. Важливу роль відіграє інтерполяція і в дослідженні багатовимірних диференціальних операторів і крайових задач у вказаних просторах \cite{Triebel80, Triebel86, LionsMagenes71}.

Хоча сфера застосувань цих класичних просторів є гігантською, в останні десятиліття виникло чимало задач, вивчення яких потребує шкал просторів, градуйованих істотно більш тонко за допомогою функціональних параметрів. Якщо останні залежать від частотних змінних, тобто характеризують гладкість (або регулярність) розподілів у термінах поведінки на нескінченності їх перетворення Фур'є, то відповідні класи розподілів називають просторами узагальненої гладкості (див., наприклад, огляди \cite{Lizorkin86, KalyabinLizorkin87} та монографії \cite[с.~53]{Triebel06} і \cite[с.~54]{Triebel10}). У випадку степеневого функціонального параметра отримуються класичні простори. Одним із перших, хто вказав на потребу у просторах узагальненої гладкості, був Л.~Хермандер \cite{Hermander65}, який увів і дослідив широкі класи таких просторів і навів їх важливі застосування до диференціальних рівнянь з частинними похідними (див. також \cite{Hermander86}). Нині різні простори узагальненої гладкості широко використовуються у теорії функцій і теорії стохастичних процесів, мають важливі застосування до багатовимірних крайових задач та інтегральних рівнянь (див. недавні монографії [9, 14, 17--23]).

Звісно, серед просторів узагальненої гладкості найбільш зручні для застосувань ті, що отримуються інтерполяцією їх класичних аналогів, при цьому природно потрібні методи, які використовують функціональні параметри інтерполяції. Ці простори успадковують важливі властивості свої класичних аналогів, зокрема, допускають локалізацію й тому їх можна  коректно означити на гладких многовидах за допомогою локальних карт.
У категорії гільбертових просторів такі простори узагальненої гладкості  виділені В.~А.~Михайлецем і О.~О.~Мурачем [22--27] та застосовані до еліптичних рівнянь та еліптичних крайових задач. Ґ.~Трібель у монографії \cite[п.~1.3.3]{Triebel10} запропонував програму дослідження і застосувань (квазі)нормованих просторів логарифмічної гладкості, де, зокрема, вказується на роль інтерполяції у їх вивченні та можливі застосування до еліптичних крайових задач.

Мета цієї статті~--- виділити і дослідити класи нормованих або квазінормованих просторів узагальненої гладкості, які можна отримати різними методами інтерполяції просторів Соболєва, Нікольського\,--\,Бєсова і Лізоркіна\,--\,Трібеля. Показником  регулярності для досліджуваних  просторів служить довільна додатна функція, RO-змінна на нескінченності за В.~Г.~Авакумовичем. Вони є широкими узагальненнями вказаних класичних просторів та допускають коректне означення на гладких многовидах й тому придатні для вивчення багатовимірних крайових задач. Використаний клас функціональних показників регулярності є, мабуть, максимально широким для вказаної мети. У статті наведено застосування узагальнених просторів Соболєва\,--\,Лізоркіна\,--\,Трібеля і Нікольського\,--\,Бєсова до еліптичних рівнянь та еліптичних крайових задач на гладких многовидах.

\medskip

\textbf{2. Класичні функціональні простори та їх інтерполяція.}
У цьому розділі розглядаємо повні нормовані функціональні
\begin{itemize}
\item[---] простори Соболєва $H^{s}_{p}$, де $s\in\mathbb{R}$ і $p\in(1,\infty)$;
\item[---] простори Нікольського\,--\,Бєсова $B^{s}_{p,q}$, де $s\in\mathbb{R}$ і $p,q\in[1,\infty]$;
\item[---] простори Лізоркіна\,--\,Трібеля $F^{s}_{p,q}$, де $s\in\mathbb{R}$ і $p,q\in[1,\infty)$.
\end{itemize}
Тут $s$~--- показник регулярності (або, інакше кажучи, гладкості) простору, а $p$ і $q$~--- показники сумовності (відповідно, основний і додатковий). (В англомовній літературі зустрічається така термінологія: “smoothness index” для $s$, “integral-exponent” для $p$ і “sum-exponent” для $q$; див., наприклад, \cite[п.~2.2, с.~35]{Johnsen96}. Простори Нікольського\,--\,Бєсова називають “Besov spaces”, а простори Лізоркіна\,--\,Трібеля~--- “Triebel--Lizorkin spaces”.) Вказані простори спочатку означаються на $\mathbb{R}^{n}$ за допомогою перетворення Фур'є \cite[п.~2.3]{Triebel80}, а потім на евклідових областях (як звуження на них просторів, заданих на $\mathbb{R}^{n}$) і компактних нескінченно гладких многовидах (за допомогою локальних карт і розбиття одиниці) \cite[п.~3.2.2]{Triebel86}. Важливо, що відповідні простори на цих многовидах не залежать від вибору локальних карт і розбиття одиниці \cite[п.~3.2.3]{Triebel86}.

Класи цих просторів пов'язані так \cite[теореми 2.3.2(d) і 2.3.3(a)]{Triebel80}:
\begin{gather}
H^{s}_{2}=B^{s}_{2,2}=F^{s}_{2,2}\quad
\mbox{для довільного}\;\;s\in\mathbb{R};\label{H-B-F}
\\ \label{H=F}
H^{s}_{p}=F^{s}_{p,2}\quad
\mbox{для усіх}\;\;s\in\mathbb{R}\;\;\mbox{і}\;\;p\in(1,\infty);
\\ \label{B=F}
B^{s}_{p,p}=F^{s}_{p,p}\quad
\mbox{для усіх}\;\;s\in\mathbb{R}\;\;\mbox{і}\;\;p\in[1,\infty).
\end{gather}
Зокрема,
\begin{equation}
H^{0}_{2}=B^{0}_{2,2}=F^{0}_{2,2}=L_{2}
\end{equation}
і
\begin{equation}\label{H=F=L}
H^{0}_{p}=F^{0}_{p,2}=L_{p}\quad\mbox{для довільного}\;\;p\in(1,\infty),
\end{equation}
де $L_{p}$~--- простір Лебега функцій, сумовних з степенем $p$. Ці рівності виконуються з точністю до еквівалентності норм. Простори у формулі \eqref{H-B-F} гільбертові з точністю до еквівалентності норм. Окремий інтерес викликає простір $B^{s}_{p,\infty}$. Він є простором Нікольського, якщо $s>0$ і $p\neq\infty$. Крім того, $B^{s}_{\infty,\infty}$ збігається (з точністю до еквівалентності норм) з простором Зігмунда $\mathcal{C}^{s}$ порядку $s$, якщо $s>0$ (див. \cite[с.~118, п.~2.5.7, теорема~(ii)]{Triebel86}). Отже, $B^{s}_{\infty,\infty}$ збігається з простором Гельдера $C^{s}$, якщо додатне число $s$ неціле (при означенні просторів Гельдера і Зігмунда базовим є простір $C(\mathbb{R}^{n})$ усіх обмежених рівномірно неперервних функцій на $\mathbb{R}^{n}$).

Надалі розглядаються комплексні лінійні простори. Вказані вище класи просторів мають важливі інтерполяційні властивості. Перед тим як їх нагадати, наведемо деякі стандартні позначення, які використовуються у теорії інтерполяції нормованих просторів і операторів на них. Нехай $E_{0}$ і $E_{1}$~--- пара (комплексних) банахових просторів, які є лінійними многовидами у деякому ґаусдорфовому лінійному топологічному просторі $T$, такими, що оператори вкладення їх у $T$ неперервні. Таку пару називають інтерполяційною, а самі банахові простори сумісними.

Нехай задано параметри $\theta\in(0,1)$ і $q\in[1,\infty]$. Через $(E_{0},E_{1})_{\theta,q}$ позначається банахів простір, отриманий дійсною інтерполяцією з параметрами $\theta$ і $q$ упорядкованої пари просторів $E_{0}$ і $E_{1}$ (різні еквівалентні методи дійсної інтерполяції, зокрема, $K$-метод Петре, викладені, наприклад, у монографіях
\cite[пп. 1.3\,--\,1.8]{Triebel80}, \cite[розд.~3]{BergLefstrem80} і \cite[п~13.3]{Agranovich13}). Через $[E_{0},E_{1}]_{\theta}$ позначається банахів простір, отриманий комплексною інтерполяцією з параметром $\theta$ упорядкованої пари просторів $E_{0}$ і $E_{1}$ (її означення див., наприклад, у щойно вказаних книгах \cite[пп. 1.9]{Triebel80}, \cite[розд.~4]{BergLefstrem80} і \cite[п~13.2b]{Agranovich13}).

Банахові простори $(E_{0},E_{1})_{\theta,q}$ і $[E_{0},E_{1}]_{\theta}$ неперервно вкладаються у топологічний простір $T$. Більше того, виконуються неперервні вкладення
\begin{gather}
E_{0}\cap E_{1}\hookrightarrow(E_{0},E_{1})_{\theta,q}\hookrightarrow E_{0}+E_{1},\\
E_{0}\cap E_{1}\hookrightarrow[E_{0},E_{1}]_{\theta}\hookrightarrow E_{0}+E_{1}.
\end{gather}
Тут $E_{0}\cap E_{1}$ і $E_{0}+E_{1}:=\{u_{0}+u_{1}:u_{0}\in E_{0},u_{1}\in E_{1}\}$~--- банахові простори, наділені відповідно нормами
\begin{gather}
\|u\|_{E_{0}\cap E_{1}}:=\max\{\|u\|_{E_{0}},\|u\|_{E_{1}}\},\\
\|u\|_{E_{0}+E_{1}}:=\inf\{\|u_{0}\|_{E_{0}}+\|u_{1}\|_{E_{1}}:
u=u_{0}+u_{1},\,u_{0}\in E_{0},\,u_{1}\in E_{1}\}.
\end{gather}

Головна властивість дійсного і комплексного інтерполяційних методів така: нехай $[E_{0},E_{1}]$ і $[Q_{0},Q_{1}]$~--- інтерполяційні пари банахових просторів і нехай $L:E_{0}+E_{1}\to Q_{0}+Q_{1}$ є лінійний оператор такий, що його звуження на простори $E_{0}$ і $E_{1}$ є обмеженими операторами $L:E_{0}\to Q_{0}$ і $L:E_{1}\to Q_{1}$; тоді звуження оператора $L$ на простори  $(E_{0},E_{1})_{\theta,q}$ і $[E_{0},E_{1}]_{\theta}$ є обмеженими операторами $L:(E_{0},E_{1})_{\theta,q}\to(Q_{0},Q_{1})_{\theta,q}$ і $L:[E_{0},E_{1}]_{\theta}\to[Q_{0},Q_{1}]_{\theta}$. Інакше кажучи, ці методи є інтерполяційними функторами. Простори $(E_{0},E_{1})_{\theta,q}$ і $[E_{0},E_{1}]_{\theta}$ називають інтерполяційними щодо пари $[E_{0},E_{1}]$ (оскільки вони задовольняють цю головну властивість, розглянуту у випадку, коли $E_{0}=Q_{0}$ і $E_{1}=Q_{1}$). Дійсний $K$-метод інтерполяції та комплексний метод інтерполяції є точними типу $\theta$, тобто виконуються такі нерівності для норм операторів, які фігурують в основній властивості:
\begin{gather}\label{interp-norm-operators-R}
\|L\|_{(E_{0},E_{1})_{\theta,q}\to(Q_{0},Q_{1})_{\theta,q}}\leq
\|L\|_{E_{0}\to Q_{0}}^{1-\theta}\cdot\|L\|_{E_{1}\to Q_{1}}^{\theta},\\
\label{interp-norm-operators-C}
\|L\|_{[E_{0},E_{1}]_{\theta}\to[Q_{0},Q_{1}]_{\theta}}\leq
\|L\|_{E_{0}\to Q_{0}}^{1-\theta}\cdot\|L\|_{E_{1}\to Q_{1}}^{\theta}.
\end{gather}
Тут, звичайно, $\|L\|_{X\to Y}$ позначає норму лінійного обмеженого оператора $L:X\to Y$, де $X$ і $Y$~--- нормовані простори.

Зауважимо, що при перестановці місцями просторів у парі результат їх дійсної або комплексної інтерполяції, взагалі кажучи, змінюється:
\begin{equation}
(E_{1},E_{0})_{\theta,q}=(E_{0},E_{1})_{1-\theta,q}\quad\mbox{і}\quad
[E_{1},E_{0}]_{\theta}=[E_{0},E_{1}]_{1-\theta}.
\end{equation}

Якщо простори $E_{0}$ і $E_{1}$ гільбертові і сепарабельні, причому $E_{1}$ неперервно вкладається в $E_{0}$ а $\Lambda$~--- породжуючий оператор для їх пари, то з точністю до еквівалентності норм
\begin{equation}
(E_{0},E_{1})_{\theta,2}=[E_{0},E_{1}]_{\theta}=
\mathrm{Dom}\,\Lambda^{\theta}
\end{equation}
--- результат квадратичної інтерполяції цих просторів зі степеневим параметром $t^{\theta}$ (див. \cite[пп. 14.1, 15]{LionsMagenes71}). Тут область визначення $\mathrm{Dom}\,\Lambda^{\theta}$ оператора $\Lambda^{\theta}$ наділена гільбертовою нормою графіка $(\|u\|_{E_{0}}^{2}+\|\Lambda^{\theta}u\|_{E_{0}}^{2})^{1/2}$, де $u\in\mathrm{Dom}\,\Lambda^{\theta}$.

Нагадаємо основні інтерполяційні властивості просторів Соболєва, Нікольського\,--\,Бєсова і Лізоркіна\,--\,Трібеля. Нехай $G$ позначає або евклідів простір $\mathbb{R}^{n}$, де $n\in\mathbb{N}$, або відкритий евклідів півпростір (зокрема, $\mathbb{R}^{n}_{+}:=\{(x',x_{n}):x'\in\mathbb{R}^{n-1},x_{n}>0\}$), або нескінченно гладкий компактний многовид без краю, або внутрішню частину
нескінченно гладкого компактного многовиду з краєм (зокрема, відкриту обмежену евклідову область класу $C^{\infty}$, див, наприклад,    \cite[п.~3.2.1, означення~2]{Triebel80}). Будемо розглядати ці функціональні простори, задані на множині $G$, позначаючи їх через $H^{s}_{p}(G)$, $B^{s}_{p,q}(G)$ і $F^{s}_{p,q}(G)$ відповідно.

Почнемо з дійсної інтерполяції цих просторів. Нехай довільно задані параметри інтерполяції $\theta\in(0,1)$ і $q\in[1,\infty]$, показники регулярності $s_{0},s_{1}\in\mathbb{R}$ та показники сумовності $p,q_{0},q_{1}\in[1,\infty]$ функціональних просторів. Коли формула містить $F$-простори, додатково припускаємо, що $p,q_{0},q_{1}\neq\infty$, а коли містить $H$-простори, додатково припускаємо, що $p\notin\{1,\infty\}$. Покладемо
\begin{equation}\label{s}
s:=(1-\theta)s_{0}+\theta s_{1}.
\end{equation}
З точністю до еквівалентності норм правильні такі рівності:
\begin{gather}\label{R-interp-B}
(B^{s_0}_{p,q_{0}}(G),B^{s_1}_{p,q_{1}}(G))_{\theta,q}=
B^{s}_{p,q}(G)\quad\mbox{при}\quad s_{0}\neq s_{1},
\\ \label{R-interp-F}
(F^{s_0}_{p,q_{0}}(G),F^{s_1}_{p,q_{1}}(G))_{\theta,q}=
B^{s}_{p,q}(G)\quad\mbox{при}\;\;s_{0}\neq s_{1}.
\end{gather}
Ці та подібні формули називають інтерполяційними. З огляду на \eqref{H=F}, окремим (і важливим) випадком останньої інтерполяційної формули є рівність
\begin{equation}\label{R-interp-H}
(H^{s_0}_{p}(G),H^{s_1}_{p}(G))_{\theta,q}=
B^{s}_{p,q}(G)\quad\mbox{при}\quad s_{0}\neq s_{1}.
\end{equation}
Отже, простори Нікольського\,--\,Бєсова отримуються дійсної інтерполяцією пар просторів Соболєва.

Перейдемо до комплексної інтерполяції розглянутих просторів. Нехай довільно задані параметр інтерполяції $\theta\in(0,1)$, показники регулярності $s_{0},s_{1}\in\mathbb{R}$ та показники сумовності $p_0,p_1,q_{0},q_{1}\in[1,\infty)$ просторів. Означимо параметр $s$ за формулою \eqref{s}, а параметри $p,q\in[1,\infty)$ за формулами
\begin{equation}\label{pq-interp}
\frac{1}{p}=\frac{1-\theta}{p_0}+\frac{\theta}{p_1}\quad\mbox{і}\quad
\frac{1}{q}=\frac{1-\theta}{q_0}+\frac{\theta}{q_1}.
\end{equation}
З точністю до еквівалентності норм правильні такі рівності:
\begin{gather}\label{C-interp-B}
[B^{s_0}_{p_{0},q_{0}}(G),B^{s_1}_{p_{1},q_{1}}(G)]_{\theta}=
B^{s}_{p,q}(G),
\\ \label{C-interp-F}
[F^{s_0}_{p_{0},q_{0}}(G),F^{s_1}_{p_{1},q_{1}}(G)]_{\theta}=
F^{s}_{p,q}(G).
\end{gather}
З огляду на \eqref{H=F}, окремим (і важливим) випадком останньої інтерполяційної формули є рівність
\begin{equation}\label{C-interp-H}
[H^{s_0}_{p_0}(G),H^{s_1}_{p_1}(G)]_{\theta}=H^{s}_{p}(G),
\end{equation}
де $p_0,p_1\neq1$. Отже, кожен з класів просторів Соболєва, Бєсова (де $p$ і $q$ скінченні) і Лізоркіна\,--\,Трібеля замкнений відносно комплексної інтерполяції. Клас просторів Нікольського не замкнений відносно комплексної інтерполяції, оскільки
\begin{equation}
[B^{s_0}_{p,\infty}(\mathbb{R}^{n}),
B^{s_1}_{p,\infty}(\mathbb{R}^{n})]_{\theta}\subsetneqq
B^{s}_{p,\infty}(\mathbb{R}^{n}),\quad\mbox{якщо}\;\;s_{0}\neq s_{1}\;\;\mbox{і}\;\;1<p<\infty
\end{equation}
(див. \cite[п.~2.4.1, теорема~(e)]{Triebel80}).

У випадку $G=\mathbb{R}^{n}$ інтерполяційні формули \eqref{R-interp-B} і \eqref{R-interp-F} доведені, наприклад, в
\cite[с.~83, п.~2.4.2, теорема]{Triebel86}, формула \eqref{C-interp-B}~--- в \cite[с.~218, п.~2.4.1, теорема~(d)]{Triebel80}, а формула \eqref{C-interp-F}~--- в \cite[с.~222, п.~2.4.2, теорема 1(d)]{Triebel80}. Правда, у цитованій монографії \cite{Triebel80} формула \eqref{C-interp-B} доведена  за припущення, що $p_0,p_1\neq1$, а формула \eqref{C-interp-F}~--- за припущення, що $p_0,p_1,q_{0},q_{1}\neq\nobreak1$. Цих припущень можна позбутися, як це зроблено в \cite[с.~195, теорема 6.4.5, формула~(7)]{BergLefstrem80} для рівності~\eqref{C-interp-B}. Для $G=\mathbb{R}^{n}_{+}$ ці формули випливають з попереднього випадку та тверджень про існування обмежених лінійних операторів продовження розподілів з просторів $B$ (Нікольського\,--\,Бєсова) і просторів $F$ (Лізоркіна\,--\,Трібеля) на $\mathbb{R}^{n}_{+}$ в аналогічні простори на $\mathbb{R}^{n}$ (порівняти, наприклад, з \cite[доведення теореми 3.2]{MikhailetsMurach10} або з \cite[доведення теореми 3.3]{06UMJ3}). Ці твердження доведені в \cite[с. 238, п.~2.9.1, твердження~2, п.~2.9.2, теорема]{Triebel86}. Зокрема, інтерполяційні формули \eqref{R-interp-B} (при $1<p<\infty$) і \eqref{C-interp-B} (при $p_0,p_1\neq1$) для $G=\mathbb{R}^{n}_{+}$ доведені в \cite[п.~2.10.1, лема]{Triebel80}. Інтерполяційні формули \eqref{R-interp-B}, \eqref{R-interp-F} і \eqref{C-interp-B}, \eqref{C-interp-F} для просторів на гладких компактних многовидах виводяться з аналогічних формул для $G\in\{\mathbb{R}^{n},\mathbb{R}^{n}_{+}\}$ за допомогою локальних карт і розбиття одиниці (порівняти з \cite[доведення теореми 2.2]{MikhailetsMurach10} у випадку многовиду без краю і з \cite[доведення теореми 3.5]{06UMJ3} у випадку многовиду з краєм). Зокрема, в \cite[п.~3.3.6, твердження~(i)]{Triebel86} доведено формули \eqref{R-interp-B} і \eqref{R-interp-F} у випадку, коли $G$~--- межа обмеженої евклідової області класу~$C^{\infty}$ (тоді $G$~--- нескінченно гладкий компактний многовид без краю).

Як бачимо, у дійсних інтерполяційних формулах \eqref{R-interp-B}\,--\,\eqref{R-interp-H} основний показник сумовності $p$ не змінюється, а $s_{0}\neq s_{1}$, на відміну від комплексних інтерполяційних формул \eqref{C-interp-B}\,--\,\eqref{C-interp-H}. Якщо змінювати його у лівій частині формул \eqref{R-interp-B}\,--\,\eqref{R-interp-H} або узяти $s_{0}=s_{1}$, отримуємо у результаті дійсної інтерполяції простори, взагалі кажучи, відмінні від просторів Соболєва, Бєсова і Лізоркіна\,--\,Трібеля  \cite[п.~2.4.1, формули (4) і~(5), п.~2.4.2, формули (3) і~(5)]{Triebel80}.

Окремий інтерес привертає інтерполяційна формула \eqref{C-interp-F} у випадку, коли $p_{0}>1$, $q_{0}=2$ і $p_{1}=q_{1}$. Тоді вона набирає вигляду
\begin{equation}
[H^{s_0}_{p_{0}}(G),B^{s_1}_{p_{1},p_{1}}(G)]_{\theta}=
F^{s}_{p,q}(G)
\end{equation}
на підставі \eqref{H=F} і \eqref{B=F}. Але на цьому шляху, інтерполюючи пари просторів Соболєва і Бєсова, можна отримати, лише ті простори Лізоркіна\,--\,Трібеля $F^{s}_{p,q}(G)$, які задовольняють умову $1/(2p)<1/q<1/2+1/(2p)$ \cite[п.~2.10.1]{Triebel80}.

\medskip

\textbf{2. Простори узагальненої гладкості.} Розглянемо узагальнення просторів Соболєва $H^{s}_{p}$, Нікольського\,--\,Бєсова $B^{s}_{p,q}$ і Лізоркіна\,--\,Трібеля $F^{s}_{p,q}$ на випадок, коли:
\begin{itemize}
  \item[---] показником регулярності служить замість числа $s$ досить загальна функція $\alpha:\nobreak[1,\infty)\to(0,\infty)$,
  \item[---] відповідні банахові простори $H^{\alpha}_{p}$, $B^{\alpha}_{p,q}$ і $F^{\alpha}_{p,q}$ отримуються у результаті комплексної або дійсної інтерполяції гільбертових просторів $H^{\beta}:=H^{\beta}_{2}$, які утворюють розширену соболєвську шкалу \cite[п.~2.4.2]{MikhailetsMurach10}, і просторів Соболєва, Нікольського\,--\,Бєсова або Лізоркіна\,--\,Трібеля відповідно.
\end{itemize}

Показник регулярності $\beta$ для простору $H^{\beta}$ пробігає деякий функціональний клас $\mathrm{RO}$, який, виявляється, можна узяти як клас показників регулярності для просторів $H^{\alpha}_{p}$, $B^{\alpha}_{p,q}$ і $F^{\alpha}_{p,q}$.

\smallskip

\textbf{Означення 1.} Клас $\mathrm{RO}$ складається з усіх вимірних за Борелем функцій $\alpha:[1,\infty)\rightarrow(0,\infty)$, для кожної з яких існують числа $b>1$ і $c\geq1$ такі, що
\begin{equation}
c^{-1}\leq\frac{\alpha(\lambda t)}{\alpha(t)}\leq c\quad
\mbox{для усіх}\;\;t\geq1\;\;\mbox{і}\;\;\lambda\in[1,b].
\end{equation}

У цьому означенні припускається, що числа $b$ і $c$ можуть залежати від $\alpha$. (Втім, якщо вважати $b$ однаковим для усіх $\alpha$, то буде, звісно, означено той самий клас $\mathrm{RO}$). Вказані в означенні функції називають RO-змінними за В.~Г.~Авакумовичем на нескінченності і є добре дослідженими (див., наприклад, \cite{Seneta85}, сс. 86--99]). Зокрема, до класу $\mathrm{RO}$ належить будь-яка неперервна функція $\alpha:[1,\infty)\rightarrow(0,\infty)$ вигляду
$$
\alpha(t):=t^{r}(\log t)^{r_{1}}(\log\log
t)^{r_{2}}\ldots(\underbrace{\log\ldots\log}_{k\;\mbox{\small разів}}
t)^{r_{k}},\quad t\gg1,
$$
де довільно вибрано ціле число $k\geq1$ і дійсні числа  $r,r_{1},\ldots,r_{k}$.

Нехай $n\in\mathbb{N}$. Слідуючи за \cite{Merucci84, CobosFernandez88}, означимо простори $B^{\alpha}_{p,q}(\mathbb{R}^{n})$ і $F^{\alpha}_{p,q}(\mathbb{R}^{n})$, де $\alpha\in\mathrm{RO}$, а показники сумовності $p$ і $q$~--- такі як у просторах Нікольського\,--\,Бєсова або Лізоркіна\,--\,Трібеля відповідно (простір $H^{\alpha}_{p}(\mathbb{R}^{n})$ розглянемо пізніше). Для цього знадобляться спеціальні послідовності функцій класу $\mathcal{S}(\mathbb{R}^{n})$. Як звичайно, $\mathcal{S}(\mathbb{R}^{n})$~--- лінійний топологічний простір Л.~Шварца усіх нескінченно диференційовних функцій $\eta:\mathbb{R}^{n}\to\mathbb{C}$, які швидко спадають на нескінченності разом з усіма частинними похідними $\partial^{\mu}\eta$, тобто $|x|^{m}\partial^{\mu}\eta(x)\to0$ при $|x|\to\infty$ для довільного числа $m>0$ і кожного мультиіндексу $\mu$ (з невід'ємними цілими компонентами). Функції класу $\mathcal{S}(\mathbb{R}^{n})$ називають тест-функціями. Нехай $\widehat{\eta}:=\mathcal{F}\eta$ позначає перетворення Фур'є функції $\eta\in\mathcal{S}(\mathbb{R}^{n})$; воно встановлює ізоморфізм простору $\mathcal{S}(\mathbb{R}^{n})$ на себе.

\smallskip

\textbf{Означення 2.} Послідовність $(\phi_j)_{j=0}^{\infty}$ тест-функцій $\phi_j\in\mathcal{S}(\mathbb{R}^{n})$ називаємо \emph{спеціальною}, якщо для деякого $N\in\mathbb{N}$ їх перетворення Фур'є задовольняють такі п'ять властивостей:
\begin{itemize}
\item[(i)] кожна функція $\widehat{\phi}_{j}\geq0$;
\item[(ii)] $\mathrm{supp}\,\widehat{\phi_{j}}\subset
    \{\xi\in\mathbb{R}^{n}:2^{j-N}\leq|\xi|\leq2^{j+N}\}$ для довільного $j\in\mathbb{N}$ та, крім того,  $\mathrm{supp}\,\widehat{\phi}_{0}\subset\{\xi\in\mathbb{R}^{n}:
|\xi|\leq2^{N}\}$;
\item[(iii)] для довільного числа $\varepsilon\in(0,1)$ існує число $c_{\varepsilon}>0$ таке, що для кожного $j\in\mathbb{N}$ виконується $\widehat{\phi}_{j}(\xi)\geq c_{\varepsilon}$, якщо $(2-\varepsilon)^{-N}2^{j}\leq|\xi|\leq(2-\varepsilon)^{N}2^{j}$, та, крім того, $\widehat{\phi}_{0}(\xi)\geq c_{\varepsilon}$, якщо $|\xi|\leq(2-\varepsilon)^{N}$.
\item[(iv)] існує число $c>0$ таке, що $c\leq\sum_{j=0}^{\infty}\widehat{\phi}_{j}(\xi)$ для довільного $\xi\in\mathbb{R}^{n}$;
\item[(v)] для кожного мультиіндексу $\mu$ існує число $\widetilde{c}_{\mu}>0$ таке, що  $2^{j|\mu|}\cdot|\partial^{\mu}\widehat{\phi}_{j}(\xi)|\leq \widetilde{c}_{\mu}$ для довільних  $\xi\in\mathbb{R}^{n}$ і $j\in\mathbb{N}$.
\end{itemize}

Приклад спеціальної послідовності, отримаємо, якщо довільно виберемо число $\nobreak{N\in\mathbb{N}}$ і функцію $\eta\in\mathcal{S}(\mathbb{R}^{n})$ таку, що $\mathrm{supp}\,\eta\subset
\{\xi\in\mathbb{R}^{n}:2^{-N}\leq|\xi|\leq2^{N}\}$ і $\eta(\xi)>0$ за умови $2^{-N}<|\xi|<2^{N}$, покладемо $\widehat{\phi}_{j}(\xi):=\eta(2^{-j}\xi)$ для усіх $\xi\in\mathbb{R}^{n}$ і $j\in\mathbb{N}$, довільно виберемо функцію $\widehat{\phi}_{0}\in\mathcal{S}(\mathbb{R}^{n})$ таку, що $\mathrm{supp}\,\widehat{\phi}_{0}\subset\{\xi\in\mathbb{R}^{n}:
|\xi|\leq2^{N}\}$ і $\widehat{\phi}_{0}>0$ за умови $|\xi|<2^{N}$ та перейдемо до Фур'є-прообразів функцій $\widehat{\phi}_{j}$ і $\widehat{\phi}_{0}$.

\smallskip

Довільно виберемо спеціальну послідовність $(\phi_j)_{j=0}^{\infty}\subset\mathcal{S}(\mathbb{R}^{n})$.

\smallskip

\textbf{Означення 3.} Нехай $\alpha\in\mathrm{RO}$, $p,q\in[1,\infty]$. Лінійний простір $B^{\alpha}_{p,q}(\mathbb{R}^{n})$ складається з усіх розподілів $f\in\mathcal{S}'(\mathbb{R}^{n})$, які задовольняють умову
\begin{equation}\label{def-B}
\|f\|_{B^{\alpha}_{p,q}(\mathbb{R}^{n})}:=
\biggl(\,\sum_{j=0}^{\infty}\alpha^{q}(2^{j})\,
\|\phi_{j}\ast f\|_{L_{p}(\mathbb{R}^{n})}^{q}\biggr)^{1/q}<\infty,
\quad\mbox{якщо}\;\;q<\infty,
\end{equation}
або умову
\begin{equation}\label{def-B-inf}
\|f\|_{B^{\alpha}_{p,q}(\mathbb{R}^{n})}:=\sup_{0\leq j\in\mathbb{Z}}
\alpha(2^{j})\,\|\phi_{j}\ast f\|_{L_{p}(\mathbb{R}^{n})}<\infty,
\quad\mbox{якщо}\;\;q=\infty.
\end{equation}
Цей простір наділено нормою $\|\cdot\|_{B^{\alpha}_{p,q}(\mathbb{R}^{n})}$. Його називаємо \emph{узагальненим простором Нікольського\,--\,Бєсова} порядку $\alpha$ (зокрема,~--- узагальненим простором Нікольського, якщо $p\neq\infty$ і $q=\infty$, або узагальненим простором Зігмунда, якщо $p=q=\infty$).

\smallskip

Тут, як звичайно, $\mathcal{S}'(\mathbb{R}^{n})$~--- лінійний топологічний простір усіх повільно зростаючих розподілів на $\mathbb{R}^{n}$ (він є дуальним до простору $\mathcal{S}(\mathbb{R}^{n})$), а $\phi_{j}\ast f$~--- згортка тест-функції $\phi_{j}\in\mathcal{S}(\mathbb{R}^{n})$ і розподілу $f\in\mathcal{S}'(\mathbb{R}^{n})$ (тому $\phi_{j}\ast f\in C^{\infty}(\mathbb{R}^{n})$). Нагадаємо, що
\begin{equation}
\|u\|_{L_{p}(\mathbb{R}^{n})}:=\|u\|_{L_{p}(\mathbb{R}^{n},dx)}:=
\biggl(\;\int\limits_{\mathbb{R}^{n}}|u(x)|^{p}dx\biggr)^{1/p},
\end{equation}
де $p<\infty$, та
\begin{equation}
\|u\|_{L_{\infty}(\mathbb{R}^{n})}:=
\|u\|_{L_{\infty}(\mathbb{R}^{n},dx)}:=
\mathrm{ess}\sup\limits_{x\in\mathbb{R}^{n}}\,|u(x)|
\end{equation}
для довільної вимірної (зокрема, неперервної) функції $u:\mathbb{R}^{n}\to\mathbb{C}$.

Норми у формулах \eqref{def-B} і \eqref{def-B-inf} можна записати одноманітно, скориставшись стандартним позначенням норми
\begin{equation}
\|\omega\|_{l_{q}}:=\left\{
\begin{array}{ll}
\biggl(\,\sum\limits_{j=0}^{\infty}|\omega_{j}|^{q}\biggr)^{1/q},&
\hbox{якщо}\;1\leq q<\infty, \\
\sup\limits_{0\leq j\in\mathbb{Z}}|\omega_{j}|,&\hbox{якщо}\;q=\infty,
\end{array}
\right.
\end{equation}
числової послідовності $\omega:=(\omega_{j})_{j=0}^{\infty}$. А саме,
\begin{equation}
\|f\|_{B^{\alpha}_{p,q}(\mathbb{R}^{n})}:=
\|(\alpha(2^{j})\,\|\phi_{j}\ast f\|_{L_{p}(\mathbb{R}^{n})})_{j=0}^{\infty}\|_{l_{q}}
\end{equation}
яке б не було $q\in[1,\infty]$. Переставивши тут місцями норми $\|\cdot\|_{L_{p}(\mathbb{R}^{n})}$ і $\|\cdot\|_{l_{q}}$, отримаємо означення простору $F^{\alpha}_{p,q}(\mathbb{R}^{n})$.

\smallskip

\textbf{Означення 4.} Нехай $\alpha\in\mathrm{RO}$, а $p,q\in[1,\infty)$. Лінійний простір $F^{\alpha}_{p,q}(\mathbb{R}^{n})$ складається з усіх розподілів $f\in\mathcal{S}'(\mathbb{R}^{n})$, які задовольняють умову
\begin{equation}\label{def-F}
\|f\|_{F^{\alpha}_{p,q}(\mathbb{R}^{n})}:=
\|\,\|(\alpha(2^{j})(\phi_{j}\ast f)(x))_{j=0}^{\infty}\|_{l_{q}}
\|_{L_{p}(\mathbb{R}^{n},dx)}<\infty.
\end{equation}
Цей простір наділено нормою $\|\cdot\|_{F^{\alpha}_{p,q}(\mathbb{R}^{n})}$. Його називаємо \emph{узагальненим простором Лізоркіна\,--\,Трібеля} порядку~$\alpha$.

\smallskip

Як бачимо,
\begin{equation}
\|f\|_{F^{\alpha}_{p,q}(\mathbb{R}^{n})}=
\biggl(\;\int\limits_{\mathbb{R}^{n}}\biggl(\,
\sum_{j=0}^{\infty}\alpha^{q}(2^{j})\,|(\phi_{j}\ast f)(x)|^{q}
\biggr)^{p/q}dx\biggr)^{1/p}.
\end{equation}

Означення~3 наведено у роботі \cite[с.~191, означення~6]{Merucci84}. Означення~4 сформульовано у більш пізній статті \cite[с.~160, означення~2.2]{CobosFernandez88} для $p,q\in(1,\infty)$ (там же наведено і означення~3). У цих двох роботах замість $\mathrm{RO}$ використано клас $\mathfrak{B}$ усіх неперервних функцій $\alpha:(0,\infty)\to(0,\infty)$ таких, що $\alpha(1)=1$ і
\begin{equation}\label{dilation-function}
\overline{\alpha}(\lambda):=
\sup_{t>0}\frac{\alpha(\lambda t)}{\alpha(t)}<\infty
\quad\mbox{для кожного}\;\;\lambda>0.
\end{equation}
Функцію $\overline{\alpha}(\lambda)$ аргументу $\lambda>0$ називають функцією розтягу функції $\alpha$. Втім, використання усіх $\alpha\in\mathrm{RO}$ або усіх $\alpha\in\mathfrak{B}$ дає той самий (з точністю до еквівалентності норм) клас просторів $B^{\alpha}_{p,q}(\mathbb{R}^{n})$ чи $F^{\alpha}_{p,q}(\mathbb{R}^{n})$. Справа у тому, що означення цих просторів використовує лише значення $\alpha(2^{j})$, де $0\leq j\in\mathbb{Z}$, а кожна з умов $\alpha\in\mathrm{RO}$ і $\alpha\in\mathfrak{B}$ тягне за собою слабку еквівалентність  $\alpha(2^{j+1})\asymp\alpha(2^{j})$ при $0\leq j\in\mathbb{Z}$, тобто властивість
\begin{equation}
(\exists c\geq1)\,(\forall j\in\mathbb{N}\cup\{0\}):\;
c^{-1}\leq\frac{\alpha(2^{j+1})}{\alpha(2^{j})}\leq c.
\end{equation}
Крім того, узявши довільну числову послідовність $\omega:=(\omega_j)_{j=0}^{\infty}\subset(0,\infty)$ таку, що $\omega_{j+1}\asymp\omega_{j}$ при $0\leq j\in\mathbb{Z}$, будуємо неперервну кусково-лінійну функцію $\alpha\in\mathrm{RO}$ таку, що $\alpha(2^{j})=\omega_{j}$ для кожного $j\in\mathbb{N}\cup\{0\}$, а за $\alpha$~--- функцію $\alpha^{\star}\in\mathfrak{B}$ за формулою
\begin{equation}\label{alpha-star}
\alpha^{\star}(t):=\left\{
\begin{array}{ll}
\alpha(t)/\alpha(1),&\hbox{якщо}\;t\geq1,\\
\alpha(1)/\alpha(t^{-1}),&\hbox{якщо}\;0<t<1.
\end{array}
\right.
\end{equation}
Отже, будь-яка вказана числова послідовність $\omega$ породжує рівні (з точністю до пропорційності норм) простори $B^{\alpha}_{p,q}(\mathbb{R}^{n})=
B^{\alpha^{\star}}_{p,q}(\mathbb{R}^{n})$ і $F^{\alpha}_{p,q}(\mathbb{R}^{n})=
F^{\alpha^{\star}}_{p,q}(\mathbb{R}^{n})$. Звісно, в означеннях цих просторів замість функціонального параметру $\alpha\in\mathrm{RO}$ можна використовувати послідовність $\omega$, як це зроблено, в \cite[означення~2.5]{HaroskeMoura08}. З огляду на зв'язок цих просторів із розширеною соболєвською шкалою $\{H^{\beta}(\mathbb{R}^{n}):\beta\in\mathrm{RO}\}$, краще використовувати функціональний параметр з класу $\mathrm{RO}$.

Простори, уведені в означеннях 3 і 4, (а також їх аналоги для більш широких класів показників регулярності $\alpha$), часто називають просторами узагальненої гладкості (spaces of generalized (or generalised) smoothness). Такий термін прийнято в оглядах \cite{Lizorkin86, KalyabinLizorkin87} робіт російських математиків за цією тематикою, в основному Г.~А.~Калябіна і М.~Л.~Гольдмана, та використовується у недавніх роботах (див., наприклад, статті \cite{HaroskeMoura08, FarkasLeopold06,MouraNevesSchneider14} і монографії Ґ.~Трібеля \cite[с.~53]{Triebel06} і \cite[с.~54]{Triebel10}). В~оглядах \cite{Lizorkin86, KalyabinLizorkin87} розглянуто простори узагальненої гладкості $B^{(\omega,\nu)}_{p,q}(\mathbb{R}^{n})$ і $F^{(\omega,\nu)}_{p,q}(\mathbb{R}^{n})$ та їх аналоги для евклідових областей, де $p,q\in(1,\infty)$, числова послідовність $\omega:=(\omega_j)_{j=0}^{\infty}\subset(0,\infty)$ задовольняє властивості
\begin{equation}
\sup_{0\leq j\in\mathbb{Z}}\frac{\omega_{j+1}}{\omega_j}<\infty\quad
\mbox{і}\quad \sum_{j=0}^{\infty}|\omega_j|^{-q'}<\infty,
\end{equation}
а числова послідовність $\nu:=(\nu_j)_{j=0}^{\infty}\subset(0,\infty)$
є сильно зростаючою у сенсі \cite[п.~Д.1.2, c.~386, означення~1]{Lizorkin86} (як звичайно, $1/q+1/q'=1$). Зокрема, якщо $\omega_{j+1}\asymp\omega_{j}$ і кожне $\nu_j=2^{j}$, то ці простори збігаються з точністю до еквівалентності норм з відповідно просторами $B^{\alpha}_{p,q}(\mathbb{R}^{n})$ і $F^{\alpha}_{p,q}(\mathbb{R}^{n})$ для деякого $\alpha\in\mathrm{RO}$ такого, що (принаймні) $\alpha(t)\to\infty$ при $t\to\infty$. Простори $B^{(\omega,\nu)}_{p,q}(\mathbb{R}^{n})$ і $F^{(\omega,\nu)}_{p,q}(\mathbb{R}^{n})$ для довільної послідовності $\omega:=(\omega_j)_{j=0}^{\infty}\subset(0,\infty)$ такої, що $\omega_{j+1}\asymp\omega_{j}$, і досить широкого класу послідовностей $\nu:=(\nu_j)_{j=0}^{\infty}\subset(0,\infty)$ уведені (і досліджені) в
\cite[Означення 3.1.2]{FarkasLeopold06} і \cite[Означення 2.5]{HaroskeMoura08} (остання робота охоплює квазінормований випадок, коли $p,q\in(0,1)$). Якщо кожне $\nu_j=2^{j}$, то класи цих просторів збігаються з класами просторів, уведених в означеннях 3 і 4 відповідно.

Зауважимо, що в означеннях 3 і 4 згортку $\phi_{j}\ast f$ можна подати у вигляді
\begin{equation}
\phi_{j}\ast f=\mathcal{F}^{-1}\mathcal{F}(\phi_{j}\ast f)=\mathcal{F}^{-1}[\widehat{\phi}_{j}\,\mathcal{F}f],
\end{equation}
де $\mathcal{F}$ і $\mathcal{F}^{-1}$~--- оператори прямого і оберненого перетворення Фур'є, задані на просторі $\mathcal{S}'(\mathbb{R}^{n})$. Отже, ці означення вводять простори $B^{\alpha}_{p,q}(\mathbb{R}^{n})$ і  $F^{\alpha}_{p,q}(\mathbb{R}^{n})$ за допомогою перетворення Фур'є. Такий підхід названо декомпозицією диференціальних властивостей розподілу $f$ \cite[с.~54]{Triebel86}.

\smallskip

\textbf{Твердження 1.} \it Лінійні нормовані простори $B^{\alpha}_{p,q}(\mathbb{R}^{n})$ і $F^{\alpha}_{p,q}(\mathbb{R}^{n})$, уведені в означеннях~$3$ і~$4$, є повними та з точністю до еквівалентності норм не залежать від вибору послідовності $(\phi_j)_{j=0}^{\infty}\subset\mathcal{S}(\mathbb{R}^{n})$, яка задовольняє означення~$2$. \rm

\smallskip

Це твердження належить \cite[с.~161]{CobosFernandez88} (правда, при $p,q\neq1$ для $F$-просторів) і доводиться подібно до випадку звичайних (тобто класичних) просторів Нікольського\,--\,Бєсова і Лізоркіна\,--\,Трібеля (див. \cite[п.~2.3.2, доведення твердження~1; п.~2.3.3, крок~4 доведення теореми]{Triebel86}).

Очевидно, що
\begin{equation}\label{B=F-gen}
B^{\alpha}_{p,p}(\mathbb{R}^{n})=F^{\alpha}_{p,p}(\mathbb{R}^{n})
\end{equation}
з рівністю норм, де $1\leq p<\infty$.

Розглянемо три важливі окремі випадки просторів $B^{\alpha}_{p,q}(\mathbb{R}^{n})$ і $F^{\alpha}_{p,q}(\mathbb{R}^{n})$.

\smallskip

\textbf{Випадок 1 (\textit{класичні простори}).} Якщо $\alpha(t)\equiv t^{s}$ для деякого $s\in\mathbb{R}$, то $B^{\alpha}_{p,q}(\mathbb{R}^{n})$ стає, за означенням, простором Нікольського\,--\,Бєсова $B^{s}_{p,q}(\mathbb{R}^{n})$, а $F^{\alpha}_{p,q}(\mathbb{R}^{n})$~--- простором Лізоркіна\,--\,Трібеля $F^{s}_{p,q}(\mathbb{R}^{n})$. (Часто використовують інше означення цих просторів \cite[п.~2.3.1, означення~1 (a) і (b)]{Triebel80}, еквівалентне даному \cite[п.~2.3.2, теорема, п. (a) і (b)]{Triebel80}.)

\smallskip

\textbf{Випадок 2 (\textit{узагальнені простори Соболєва}).}
Припустимо, що функціональний параметр $\alpha\in\mathrm{RO}$ задовольняє такі дві умови: $\alpha\in C^{\infty}([1,\infty))$ і
\begin{equation}\label{cond-for-H-gen}
(\forall\,m\in\mathbb{N})\,(\exists\,c_{m}>0)\,
(\forall\,t\in[1,\infty)):\;t^{m}|\alpha^{(m)}(t)|\leq c_{m}\alpha(t).
\end{equation}
Нехай $p\in(1,\infty)$.
Тоді
\begin{equation}\label{def-H-alpha}
F^{\alpha}_{p,2}(\mathbb{R}^{n})=\bigl\{f\in\mathcal{S}'(\mathbb{R}^{n}):
\mathcal{F}^{-1}[\alpha(\langle\xi\rangle)(\mathcal{F}f)(\xi)]\in L_{p}(\mathbb{R}^{n})\bigr\}=:H^{\alpha}_{p}(\mathbb{R}^{n})
\end{equation}
і виконується еквівалентність норм
\begin{equation}\label{def-H-alpha-norm}
\|f\|_{F^{\alpha}_{p,2}(\mathbb{R}^{n})}\asymp
\|\mathcal{F}^{-1}[\alpha(\langle\xi\rangle)(\mathcal{F}f)(\xi)]\,\|_
{L_{p}(\mathbb{R}^{n})}=:\|f\|_{H^{\alpha}_{p}(\mathbb{R}^{n})};
\end{equation}
тут $\langle\xi\rangle:=(1+|\xi|^{2})^{1/2}$~--- згладжений модуль вектора $\xi\in\mathbb{R}^{n}$. Це доведено в \cite[с.~162, теорема 3.4]{CobosFernandez88}. Якщо $\alpha(t)\equiv t^{s}$ для деякого $s\in\mathbb{R}$, то $H^{\alpha}_{p}(\mathbb{R}^{n})$ стає простором Соболєва $H^{s}_{p}(\mathbb{R}^{n})$. Тому банахів простір $H^{\alpha}_{p}(\mathbb{R}^{n})$ називаємо узагальненим простором Соболєва порядку $\alpha$ з показником сумовності $p$.

Означення \eqref{def-H-alpha} простору  $H^{\alpha}_{p}(\mathbb{R}^{n})$ коректне, якщо функція $\alpha(\langle\xi\rangle)$ аргументу $\xi\in\mathbb{R}^{n}$ є поточковим мультиплікатором на просторі $\mathcal{S}'(\mathbb{R}^{n})$. Включення $\alpha\in\mathrm{RO}\cap C^{\infty}([1,\infty))$ разом з умовою \eqref{cond-for-H-gen} гарантує це. Доведення формули $F^{\alpha}_{p,2}(\mathbb{R}^{n})=H^{\alpha}_{p}(\mathbb{R}^{n})$, наведене в \cite[с.~162, теорема~3.4]{CobosFernandez88}, використовує оператор
\begin{equation}
J^{\alpha}:f\mapsto
\mathcal{F}^{-1}[\alpha(\langle\xi\rangle)(\mathcal{F}f)(\xi)],
\quad\mbox{де}\;\;f\in\mathcal{S}'(\mathbb{R}^{n}),
\end{equation}
рівність \eqref{H=F=L} та ізоморфізми $J^{\alpha}:H^{\alpha}_{p}(\mathbb{R}^{n})\leftrightarrow H^{0}_{p}(\mathbb{R}^{n})$ (очевидний) і $J^{\alpha}:F^{\alpha}_{p,2}(\mathbb{R}^{n})\leftrightarrow F^{0}_{p,2}(\mathbb{R}^{n})$. Останній обґрунтовується подібно до \cite[п.~2.3.4, доведення]{Triebel80}, де властивість \eqref{cond-for-H-gen} потрібна для побудови за $\alpha$ деякої спеціальної послідовності тест-функцій (див. означення~2).

Узагальнені простори Соболєва досліджувалися у роботах \cite{VolevichPaneah65, Schechter67, Triebel77III, Triebel77IV, Merucci84, CobosFernandez88} 60--80-х років минулого сторіччя, як ізотропні \cite{Merucci84, CobosFernandez88} так і, взагалі кажучи, анізотропні простори \cite{VolevichPaneah65, Schechter67, Triebel77III, Triebel77IV}, коли показником регулярності є нерадіальна функція $n$ змінних. У вказаних роботах, за винятком \cite{VolevichPaneah65} і \cite{Triebel77IV}, вивчали й інтерполяційні властивості цих просторів. На межі минулого і цього століття з'явилася серія робіт А.~Г.~Багдасаряна \cite{Bagdasaryan92, Bagdasaryan96, Bagdasaryan97, Bagdasaryan98, Bagdasaryan01, Bagdasaryan04, Bagdasaryan05, Bagdasaryan10}, у яких досліджуються інтерполяційні та інші властивості анізотропних просторів Соболєва (а також Бєсова і Лізоркіна\,--\,Трібеля) узагальненої гладкості. У новітній статті \cite{Faierman20} дано застосування (мабуть, уперше) анізотропних узагальнених просторів Соболєва  до еліптичних рівнянь на~$\mathbb{R}^{n}$.

\smallskip

\textbf{Випадок 3 (\textit{гільбертові простори}).} Їх отримуємо, коли $p=q=2$. Оскільки перетворення Фур'є є ізоморфізмом на просторі $L_{2}(\mathbb{R}^{n})$, то означення узагальненого простору Соболєва $H^{\alpha}(\mathbb{R}^{n}):=H^{\alpha}_{2}(\mathbb{R}^{n})$ і гільбертової норми у ньому можна дати для довільного $\alpha\in\mathrm{RO}$ за формулами \eqref{def-H-alpha} і \eqref{def-H-alpha-norm}, прибравши у них обернене перетворення Фур'є $\mathcal{F}^{-1}$ і вимагаючи, щоб $\mathcal{F}f$ була класичною функцією.

\smallskip

\textbf{Теорема 1.} \it Для кожного $\alpha\in\mathrm{RO}$ виконуються рівності просторів
\begin{equation}
F^{\alpha}_{2,2}(\mathbb{R}^{n})=B^{\alpha}_{2,2}(\mathbb{R}^{n})=
H^{\alpha}_{2}(\mathbb{R}^{n})
\end{equation}
з точністю до еквівалентності норм. \rm

\smallskip

Ця теорема доводиться (як і рівність \eqref{def-H-alpha}) за допомогою оператора $J^{\alpha}$, який тепер коректно означений для усіх
\begin{equation}
f\in H^{-\infty}_{2}(\mathbb{R}^{n}):=
\bigcup_{s\in\mathbb{R}}H^{s}_{2}(\mathbb{R}^{n})=
\bigcup_{\alpha\in\mathrm{RO}}H^{\alpha}_{2}(\mathbb{R}^{n})=
\bigcup_{\alpha\in\mathrm{RO}}F^{\alpha}_{2,2}(\mathbb{R}^{n}).
\end{equation}
При обґрунтуванні ізоморфізму $J^{\alpha}:F^{\alpha}_{2,2}(\mathbb{R}^{n})\leftrightarrow L_{2}(\mathbb{R}^{n})$ нема потреби у переході від стартової спеціальної системи тест-функцій до іншої такої системи, пов'язаної з $\alpha$, й тому не потрібна властивість \eqref{cond-for-H-gen}.

Клас усіх гільбертових просторів $H^{\alpha}(\mathbb{R}^{n})$, де $\alpha\in\mathrm{RO}$, виділений і досліджений в \cite{MikhailetsMurach10, MikhailetsMurach13UMJ3, MikhailetsMurach15ResMath1} та названий розширеною соболєвською шкалою. Вона отримується квадратичною інтерполяцією (з функціональним параметром) пар гільбертових просторів Соболєва та є замкненою відносно квадратичної інтерполяції. За допомогою цієї властивості для вказаної шкали побудована теорія еліптичних систем на многовидах і теорія еліптичних крайових задач \cite{ZinchenkoMurach12UMJ11, AnopMurach14UMJ, AnopKasirenko16MFAT, AnopDenkMurach20arxiv}.

\medskip

\textbf{3. Інтерполяція просторів узагальненої гладкості.} Розглянемо комплексну та дійсну інтерполяцію просторів $B^{\alpha}_{p,q}(\mathbb{R}^{n})$ і $F^{\alpha}_{p,q}(\mathbb{R}^{n})$. Почнемо з комплексної інтерполяції.

\smallskip

\textbf{Теорема 2.} \it Нехай довільно задані параметр інтерполяції $\theta\in(0,1)$, показники регулярності $\alpha_{0},\alpha_{1}\in\mathrm{RO}$ та показники сумовності $p_0,p_1,q_{0},q_{1}\in[1,\infty)$. Означимо функціональний параметр $\alpha\in\mathrm{RO}$ за формулою
\begin{equation}\label{alpha-interp}
\alpha(t):=\alpha_{0}^{1-\theta}(t)\alpha_{1}^{\theta}(t)
\quad\mbox{при}\;\;t\geq1,
\end{equation}
а параметри $p,q\in[1,\infty)$ згідно з \eqref{pq-interp}. Тоді
\begin{gather}\label{C-interp-B-gen}
[B^{\alpha_0}_{p_{0},q_{0}}(\mathbb{R}^{n}),
B^{\alpha_1}_{p_{1},q_{1}}(\mathbb{R}^{n})]_{\theta}=
B^{\alpha}_{p,q}(\mathbb{R}^{n}),\\
\label{C-interp-F-gen}
[F^{\alpha_0}_{p_{0},q_{0}}(\mathbb{R}^{n}),
F^{\alpha_1}_{p_{1},q_{1}}(\mathbb{R}^{n})]_{\theta}=
F^{\alpha}_{p,q}(\mathbb{R}^{n})
\end{gather}
з точністю до еквівалентності норм. \rm

\smallskip

З формули \eqref{C-interp-F-gen} випливає такий результат для узагальнених просторів Соболєва:

\smallskip

\textbf{Наслідок 1.} \it Додатково до умови теореми~$2$ припустимо, що функціональні параметри $\alpha_{0}$ і $\alpha_{1}$ належать до $C^{\infty}([1,\infty))$ і задовольняють умову \eqref{cond-for-H-gen}, де замість $\alpha$ узято $\alpha_{0}$ або $\alpha_{1}$. Крім того, припустимо, що $p_0,p_1\in(1,\infty)$. Тоді
\begin{equation}\label{C-interp-H-gen}
[H^{\alpha_0}_{p_{0}}(\mathbb{R}^{n}),
H^{\alpha_1}_{p_{1}}(\mathbb{R}^{n})]_{\theta}=
H^{\alpha}_{p}(\mathbb{R}^{n})
\end{equation}
з точністю до еквівалентності норм. \rm

\smallskip

Класичні інтерполяційні формули \eqref{C-interp-B}\,--\,\eqref{C-interp-H} для $G:=\mathbb{R}^{n}$ є окремими випадками формул \eqref{C-interp-B-gen}\,--\,\eqref{C-interp-H-gen} відповідно.

Зауважимо, що у статті Ґ.~Трібеля \cite[с. 239, теорема 4.2/2]{Triebel77III} доведено рівності, які мають вигляд, цілком аналогічний формулам \eqref{C-interp-B-gen} і~\eqref{C-interp-F-gen}. Ці рівності стосуються як ізотропних так і анізотропних просторів (для ізотропних просторів вимоги Ґ.~Трібеля щодо показника регулярності є іншими, ніж наша вимога належати до класу $\mathrm{RO}$). При цьому стосовно показників сумовності припускається, що $p_0,p_1\in(1,\infty)$ і $q_0,q_1\in[1,\infty)$ для $B$-просторів, та $p_0,p_1,q_0,q_1\in(1,\infty)$ для $F$-просторів. Проте, у цих рівностях $B$- і $F$-простори узагальненої гладкості, взагалі кажучи, відрізняються від просторів, уведених в означеннях~3 і~4, навіть у класичному випадку, коли $\alpha(t)$~--- степенева функція. Це показано в \cite[с. 246, теорема 6.2/1]{Triebel77III} для простору $B^{\alpha}_{2,q}(\mathbb{R}^{n})$, де $q\neq2$. Втім, згідно з \cite[с. 240, теорема 5.1/2]{Triebel77III}, простори $H^{\alpha}_{p}(\mathbb{R}^{n})$, розглянуті нами, є окремим випадком $F$-просторів, уведених у цитованій статті Ґ.~Трібеля. Отже, формула \eqref{C-interp-H-gen} є окремим випадком його результата \cite[с. 239, теорема 4.2/2]{Triebel77III}. При $p_{0}=p_{1}$ вона міститься у більш ранньому результаті М.~Шехтера \cite[с.~128, теорема~4.5]{Schechter67}.
У випадку, коли $\alpha_{0}(t)$ і $\alpha_{1}(t)$ є степенями деякої однієї функції, формула \eqref{C-interp-B-gen} міститься у недавньому результаті А.~Г.~Багдасаряна \cite[с.~45, теорема 2.4~(i)]{Bagdasaryan10}, який стосується анізотропних просторів. Якщо ця функція належить до класу $\mathrm{RO}$, то формули \eqref{C-interp-B-gen} і \eqref{C-interp-F-gen} доведені у більш ранній статті В.~Р.~Кноповой \cite[с.~641, лема~1, с.~642, зауваження~3]{Knopova06}, де розглянуті ізотропні простори узагальненої гладкості.

\smallskip

Обговоримо доведення теореми 2. Обґрунтуємо, наприклад, формулу~\eqref{C-interp-F-gen}. Нехай спочатку параметри $\alpha\in\mathrm{RO}$ і числа $p,q\in[1,\infty)$ довільні. Слідуючи за \cite[с.~160, формула~(1)]{CobosFernandez88}, розглянемо лінійний простір
\begin{gather}
l_{q}^{\alpha}:=\biggl\{\omega:=(\omega_j)_{j=0}^{\infty}\subset\mathbb{C}:
\|\omega\|_{l_{q}^{\alpha}}:=\biggl(\,\sum_{j=0}^{\infty}
\alpha^{q}(2^j)\,|\omega_j|^{q}\biggr)^{1/q}<\infty\biggr\},
\end{gather}
наділений нормою $\|\cdot\|_{l_{q}^{\alpha}}$; він є банаховим простором. Згідно з \cite[с.~161, теорема 2.5~(i)]{CobosFernandez88}, простір $F^{\alpha}_{p,q}(\mathbb{R}^{n})$ є ретрактом банахового простору $L_{p}(\mathbb{R}^{n},l_{q}^{\alpha})$, причому відповідні оператор ретракції $\mathcal{R}$ і коретракції $\mathcal{T}$ не залежать (як відображення чи правила) від $\alpha$, $p$ і~$q$. (Правда, цей результат доведено в \cite{CobosFernandez88} для $p,q\in(1,\infty)$, але у випадку, коли $p=1$ та/або $q=1$, він обґрунтовується аналогічно.) Як звичайно, якщо $E$~--- банахів простір (зокрема, $E:=l_{q}^{\alpha}$), то $L_{p}(\mathbb{R}^{n},E)$ позначає банахів простір усіх сильно вимірних відносно міри Лебега (абстрактних) функцій $g:\mathbb{R}^{n}\to E$ (точніше, класів еквівалентних функцій) таких, що
\begin{equation}
\|g\|_{L_{p}(\mathbb{R}^{n},E)}:=\biggl(\;\int\limits_{\mathbb{R}^{n}}
\|g(x)\|_{E}^{p}\,dx\biggr)^{1/p}<\infty;
\end{equation}
він наділений нормою $\|\cdot\|_{L_{p}(\mathbb{R}^{n},E)}$.

Поняття ретракту, ретракції та коретракції узято з теорії категорій й широко застосовується у теорії інтерполяції просторів (див., наприклад, монографії \cite[п.~1.2.4]{Triebel80} і \cite[п.~6.4]{BergLefstrem80}). У нашій ситуації це значить таке: існують обмежені лінійні оператори
\begin{gather}
\mathcal{R}:L_{p}(\mathbb{R}^{n},l_{q}^{\alpha})\to F^{\alpha}_{p,q}(\mathbb{R}^{n}),\\
\mathcal{T}:F^{\alpha}_{p,q}(\mathbb{R}^{n})\to L_{p}(\mathbb{R}^{n},l_{q}^{\alpha})
\end{gather}
такі, що $\mathcal{R}\mathcal{T}$ є тотожним оператором на $F^{\alpha}_{p,q}(\mathbb{R}^{n})$. (Звідси випливає, що оператор $\mathcal{R}$ сюр'єктивний, а $\mathcal{T}$ ін'єктивний.) Вони будуються подібно до класичного випадку, коли $\alpha(t)$ є степеневою функцією, причому їх побудова не залежать від параметрів $\alpha$, $p$ і~$q$; див., наприклад, \cite[п.~6.4, с.~192]{BergLefstrem80} у класичному випадку. Інтерполяція пар ретрактів просторів зводиться до інтерполяції самих просторів. Остання досліджена в \cite[пп. 1.18.4, 1.18.5]{Triebel80}. Наведемо відповідні міркування.

Нехай тепер параметри $\alpha_{0},\alpha_{1},\alpha$ і  $p_0,p_1,p,q_{0},q_{1},q$ задовольняють умову теореми~2. Оскільки простори $F^{\alpha_0}_{p_0,q_0}(\mathbb{R}^{n})$ і $F^{\alpha_1}_{p_1,q_1}(\mathbb{R}^{n})$ є ретрактами відповідно просторів $L_{p_0}(\mathbb{R}^{n},l_{q_0}^{\alpha_0})$ і $L_{p_1}(\mathbb{R}^{n},l_{q_1}^{\alpha_1})$, причому ретракція і коретракція не залежать від цих параметрів, то
\begin{equation}\label{C-interp-retract}
[F^{\alpha_0}_{p_0,q_0}(\mathbb{R}^{n}),
F^{\alpha_1}_{p_1,q_1}(\mathbb{R}^{n})]_{\theta}
\quad\mbox{є ретрактом простору}\quad
[L_{p_0}(\mathbb{R}^{n},l_{q_0}^{\alpha_0}),
L_{p_1}(\mathbb{R}^{n},l_{q_1}^{\alpha_1})]_{\theta}
\end{equation}
(з тими же самими ретракцією і коретракцією). На підставі \cite[п 1.18.4, с.~149, теорема]{Triebel80} і \cite[п 1.18.5, с.~151, теорема]{Triebel80} виконуються такі рівності просторів разом з еквівалентністю норм у них:
\begin{equation}\label{Lp-interpolation}
\bigl[L_{p_0}(\mathbb{R}^{n},l_{q_0}^{\alpha_0}),
L_{p_1}(\mathbb{R}^{n},l_{q_1}^{\alpha_1})\bigr]_{\theta}=
L_{p}\bigl(\mathbb{R}^{n},
[l_{q_0}^{\alpha_0},l_{q_1}^{\alpha_1}]_{\theta}\bigr)=
L_{p}(\mathbb{R}^{n},l_{q}^{\alpha}).
\end{equation}
(Зауважимо, що в позначеннях \cite[п 1.18.5, с.~151]{Triebel80} простір $l_{q}^{\alpha}$ є ваговим простором $L_{q,w^{q}}$, заданим на просторі $\mathbb{N}\cup\{0\}$ з дискретною мірою, де $w(j):=\alpha(2^{j})$ для кожного $j\in\mathbb{N}\cup\{0\}$, причому у цих позначеннях банахів простір $A=\mathbb{C}$.) Отже, згідно з \eqref{C-interp-retract} і \eqref{Lp-interpolation}, простір $[F^{\alpha_0}_{p_0,q_0}(\mathbb{R}^{n}),
F^{\alpha_1}_{p_1,q_1}(\mathbb{R}^{n})]_{\theta}$ є ретрактом простору $L_{p}(\mathbb{R}^{n},l_{q}^{\alpha})$. Як зазначено вище, $F^{\alpha}_{p,q}(\mathbb{R}^{n})$ є також ретрактом простору $L_{p}(\mathbb{R}^{n},l_{q}^{\alpha})$. Тому  рівність \eqref{C-interp-F-gen} правильна з точністю до еквівалентності норм.

Інша формула \eqref{C-interp-B-gen} доводиться за допомогою подібних міркувань. Порівняно з щойно наведеними міркуваннями слід замінити $L_{p}(\mathbb{R}^{n},l_{q}^{\alpha})$ на  $l_{q}^{\alpha}(L_{p}(\mathbb{R}^{n}))$. Тут, якщо $E$~--- банахів простір (зокрема, $E=L_{p}(\mathbb{R}^{n})$), то
\begin{equation}
l_{q}^{\alpha}(E):=\biggl\{\omega:=(\omega_{j})_{j=0}^{\infty}\subset E:
\|\omega\|_{l_{q}^{\alpha}(E)}:=\biggl(\,\sum_{j=0}^{\infty}
\alpha^{q}(2^j)\,\|\omega_j\|_{E}^{q}\biggr)^{1/q}<\infty\biggr\}
\end{equation}
--- банахів простір, наділений нормою $\|\cdot\|_{l_{q}^{\alpha}(E)}$. Замість \eqref{Lp-interpolation} використовується формула
\begin{equation}
\bigl[l_{q_0}^{\alpha_0}(L_{p_0}(\mathbb{R}^{n})),
l_{q_1}^{\alpha_1}(L_{p_1}(\mathbb{R}^{n}))\bigr]_{\theta}=
l_{q}^{\alpha}\bigl
([L_{p_0}(\mathbb{R}^{n}),L_{p_1}(\mathbb{R}^{n})]_{\theta}\bigr)=
l_{q}^{\alpha}(L_{p}(\mathbb{R}^{n})),
\end{equation}
правильна на підставі \cite[п 1.18.1, сc.~139 і 140, теорема]{Triebel80} і \cite[п 1.18.4, с.~149, теорема]{Triebel80} відповідно.

Для нас ключову роль відіграватиме той факт, що усі простори $B^{\alpha}_{p,q}(\mathbb{R}^{n})$ і $F^{\alpha}_{p,q}(\mathbb{R}^{n})$, де $\alpha\in\mathrm{RO}$ і $p,q\in(1,\infty)$, отримується за інтерполяційними формулами \eqref{C-interp-B-gen} і \eqref{C-interp-F-gen} відповідно, де $p_0=q_0=2$, а функція $\alpha_1(t)$ степенева. Тоді $B^{\alpha_0}_{p_0,q_0}(\mathbb{R}^{n})$ і $F^{\alpha_0}_{p_0,q_0}(\mathbb{R}^{n})$ є гільбертовим узагальненим  простором Соболєва $H^{\alpha_0}_{2}(\mathbb{R}^{n})$, а $B^{\alpha_1}_{p_1,q_1}(\mathbb{R}^{n})$ і $F^{\alpha_1}_{p_1,q_1}(\mathbb{R}^{n})$ є класичними просторами Бєсова і Лізоркіна\,--\,Трібеля.

\smallskip

\textbf{Теорема 3.} \it Нехай $\alpha\in\mathrm{RO}$ і $p,q\in(1,\infty)$. У випадку $p\neq2$ виберемо число $p_{1}$ таке, що $1<p_{1}<p<2$ або
$2<p<p_{1}<\infty$, означимо число $\theta\in(0,1)$ за формулою
\begin{equation}\label{p-2-interp}
\frac{1}{p}=\frac{1-\theta}{2}+\frac{\theta}{p_1},
\end{equation}
а число $q_1$~--- за формулою
\begin{equation}\label{q-2-interp}
\frac{1}{q}=\frac{1-\theta}{2}+\frac{\theta}{q_1}
\end{equation}
та припустимо, що $q_1\in(1,\infty)$ (це припущення виконується, зокрема, якщо $|p-p_1|$ достатньо мале). У випадку, коли $p=2$, а $q\neq2$, покладемо $p_{1}=2$, виберемо число $q_{1}$ таке, що $1<q_{1}<q<2$ або
$2<q<q_{1}<\infty$ і означимо число $\theta\in(0,1)$ за формулою~\eqref{q-2-interp}. В обох випадках довільно виберемо число $s_1\in\mathbb{R}$ і означимо функцію $\alpha_0\in\mathrm{RO}$ за формулою \begin{equation}\label{alpha-0-interp}
\alpha_{0}(t):=(t^{-\theta s_1}\alpha(t))^{1/(1-\theta)}
\quad\mbox{при}\;\;t\geq1.
\end{equation}
Тоді
\begin{gather}\label{B-gen-C-interp}
B^{\alpha}_{p,q}(\mathbb{R}^{n})=[H^{\alpha_0}_{2}(\mathbb{R}^{n}),
B^{s_1}_{p_{1},q_{1}}(\mathbb{R}^{n})]_{\theta},\\
\label{F-gen-C-interp}
F^{\alpha}_{p,q}(\mathbb{R}^{n})=[H^{\alpha_0}_{2}(\mathbb{R}^{n}),
F^{s_1}_{p_{1},q_{1}}(\mathbb{R}^{n})]_{\theta}
\end{gather}
з точністю до еквівалентності норм. \rm

\smallskip

Цей результат є наслідком теореми~2.

За формулою \eqref{B-gen-C-interp} не можна отримати простір $B^{\alpha}_{p,q}(\mathbb{R}^{n})$ у випадку, коли принаймні один з параметрів $p$ чи $q$ дорівнює $1$ або $\infty$, навіть, якщо допустити, що параметри $p_1$ чи $q_1$ можуть набувати значення $1$ або $\infty$. Так само, за цією формулою не можна отримати простір $F^{\alpha}_{p,q}(\mathbb{R}^{n})$ у випадку, коли $p=1$ та/або $q=1$. Ці випадки можна охопити, скориставшись дійсною інтерполяцією просторів, що фігурують у правій частин цієї формули. Ця інтерполяція дає узагальнені простори Нікольського\,--\,Бєсова, причому за додаткових умов на показники регулярності $\alpha_0,\alpha_1\in\mathrm{RO}$ (порівняти з \eqref{R-interp-B}).

Розглянемо версії теорем 2 і 3 для методу дійсної інтерполяції з функціональним параметром $\gamma\in\mathfrak{B}$ (клас $\mathfrak{B}$ означений там, де наведена формула \eqref{dilation-function}) і числовим параметром $q\in[1,\infty]$. Якщо $E_0$ і $E_1$ є інтерполяційна пара банахових просторів, то через $(E_0,E_1)_{\gamma,q}$ позначаємо банахів простір, отриманий вказаним методом, застосованим до цієї пари (див., означення в \cite[с.~164]{CobosFernandez88} або \cite[с.~294]{Gustavsson78}, наприклад). Вказаний метод є інтерполяційним функтором (тобто зберігає обмеженість лінійних операторів), якщо
\begin{equation}\label{Boyd-indexes}
\nu_{0}(\overline{\gamma}):=
\lim_{t\to0+}\frac{\log\overline{\gamma}(t)}{\log t}>0\quad\mbox{і}\quad
\nu_{1}(\overline{\gamma}):=
\lim_{t\to\infty}\frac{\log\overline{\gamma}(t)}{\log t}<1,
\end{equation}
де $\overline{\gamma}$~--- функція розтягу функції $\gamma$ (див. формулу \eqref{dilation-function}). За цієї умови
\begin{equation}\label{R-interp-property}
\|L\|_{(E_{0},E_{1})_{\gamma,q}\to(Q_{0},Q_{1})_{\gamma,q}}\leq
c\cdot\|L\|_{E_{0}\to Q_{0}}\cdot\overline{\gamma}
\biggl(\frac{\|L\|_{E_{1}\to Q_{1}}}{\|L\|_{E_{0}\to Q_{0}}}\biggr)
\end{equation}
у позначеннях, використаних в оцінках \eqref{interp-norm-operators-R} і
\eqref{interp-norm-operators-C}. Тут число $c>0$ не залежить від оператора $L$ та пар просторів $E_{0},E_{1}$ і $Q_{0},Q_{1}$. Це доведено в \cite[с.~295, Теорема~2.1]{Gustavsson78}; зауважимо, що умова \eqref{Boyd-indexes} еквівалентна умові \cite[с.~290, формула~(4)]{Gustavsson78}, де $f:=\gamma$, з огляду на \cite[с.~159, формули (3), (4)]{CobosFernandez88}, а умова зростання функції $\gamma$, використана в \cite[с.~290, формула~(2)]{Gustavsson78} випливає з умови $\nu_{0}(\overline{\gamma})>0$ з точністю до слабкої еквівалентності функцій на $(0,\infty)$ (якщо $\gamma$ зростає, то $c=1$). Числа $\nu_{0}(\overline{\gamma})$ і $\nu_{1}(\overline{\gamma})$ називають відповідно нижнім і верхнім індексами Бойда функції $\overline{\gamma}$ \cite{Boyd67}. У випадку, коли $\gamma(t)\equiv t^{\theta}$ для деякого $\theta\in(0,1)$, вказаний метод стає класичним дійсним методом інтерполяції з числовими параметрами $\theta$ і $q$ (у цьому випадку простори $(E_0,E_1)_{\gamma,q}$ і $(E_0,E_1)_{\theta,q}$ збігаються з рівністю норм у них).

Формула \eqref{R-interp-B} для класичних просторів правильна за умови $s_0\neq s_1$. Нам потрібен аналог цієї умови для функціональних параметрів $\alpha_0,\alpha_1\in\mathrm{RO}$, узятих замість числових параметрів $s_0$ і~$s_1$. Сформулюємо її за допомогою індексів Матушевської функцій класу $\mathrm{RO}$. Нагадаємо означення цих індексів.

Відомо \cite[додаток~1, теорема~2]{Seneta85}, що для
кожної функції $\alpha\in\mathrm{RO}$ існують дійсні числа $s_{0}<s_{1}$ і $c_{0},c_{1}>0$ такі, що
\begin{equation}\label{RO-property}
c_{0}\lambda^{s_{0}}\leq\frac{\alpha(\lambda t)}{\alpha(t)}\leq
c_{1}\lambda^{s_{1}} \quad\mbox{для всіх}\quad t,\lambda\in[1,\infty).
\end{equation}
Покладемо
\begin{gather*}
\sigma_{0}(\alpha):=
\sup\,\{s_{0}\in\mathbb{R}:\,\mbox{виконується ліва нерівність в \eqref{RO-property}}\},\\
\sigma_{1}(\alpha):=\inf\,\{s_{1}\in\mathbb{R}:\,\mbox{виконується права нерівність
в \eqref{RO-property}}\}.
\end{gather*}
Числа
$\sigma_{0}(\alpha)$ і $\sigma_{1}(\alpha)$ є відповідно нижнім і верхнім індексами Матушевської \cite{Matuszewska64} функції $\alpha\in\mathrm{RO}$ (див. також монографію \cite[п.~2.1.2]{BinghamGoldieTeugels89}). Звісно, $-\infty<\sigma_{0}(\alpha)\leq\sigma_{1}(\alpha)<\infty$. Відмітимо \cite[сс. 75, 76]{KreinPetuninSemenov78}, що
\begin{equation}\label{Matuszewska-Boyd}
\sigma_{0}(\alpha)=\nu_{0}(\overline{\alpha^{\star}})\quad\mbox{і}\quad
\sigma_{1}(\alpha)=\nu_{1}(\overline{\alpha^{\star}}),
\end{equation}
де $\overline{\alpha^{\star}}$~--- функція розтягу функції $\alpha^{\star}$, означеної за формулою \eqref{alpha-star} (якщо $\alpha\in\mathrm{RO}$ неперервна, то $\alpha^{\star}\in\mathfrak{B}$).

\smallskip

\textbf{Теорема 4.} \it Нехай задано параметри інтерполяції $\gamma\in\mathfrak{B}$ і $q\in[1,\infty]$, показники регулярності  $\alpha_{0},\alpha_{1}\in\mathrm{RO}$ та показники сумовності  $p,q_{0},q_{1}\in[1,\infty]$. Припустимо, що $\gamma$ задовольняє умову \eqref{Boyd-indexes}, а показники регулярності~--- умову
\begin{equation}\label{R-inerp-cond-index}
\sigma_{0}(\alpha_{0}/\alpha_{1})>0\quad\mbox{або}\quad
\sigma_{1}(\alpha_{0}/\alpha_{1})<0.
\end{equation}
Означимо параметр $\alpha\in\mathrm{RO}$ за формулою
\begin{equation}\label{alpha-R-interp-gamma}
\alpha(t):=\frac{\alpha_{0}(t)}{\gamma(\alpha_{0}(t)/\alpha_{1}(t))}
\quad\mbox{при}\;\;t\geq1.
\end{equation}
Тоді
\begin{equation}\label{R-interp-B-gen}
(B^{\alpha_0}_{p,q_{0}}(\mathbb{R}^{n}),
B^{\alpha_1}_{p,q_{1}}(\mathbb{R}^{n}))_{\gamma,q}=
B^{\alpha}_{p,q}(\mathbb{R}^{n}).
\end{equation}
Крім того, якщо $p,q_{0},q_{1}\in[1,\infty)$, то
\begin{equation}\label{R-interp-F-gen}
(F^{\alpha_0}_{p,q_{0}}(\mathbb{R}^{n}),
F^{\alpha_1}_{p,q_{1}}(\mathbb{R}^{n}))_{\gamma,q}=
B^{\alpha}_{p,q}(\mathbb{R}^{n}).
\end{equation}
Ці рівності просторів виконуються з точністю до еквівалентності норм у них. \rm

\smallskip

У цій теоремі умова \eqref{R-inerp-cond-index} грає роль умови $s_0\neq s_1$ у формулах \eqref{R-interp-B} і~\eqref{R-interp-F}. Звісно, якщо $\alpha_0(t)\equiv t^{s_0}$ і $\alpha_1(t)\equiv t^{s_1}$, то $s_0\neq s_1$ еквівалентно \eqref{R-inerp-cond-index}.

Якщо параметр інтерполяції $\gamma(t)\equiv t^{\theta}$ для деякого $\theta\in(0,1)$, то формула \eqref{alpha-R-interp-gamma} набирає вигляду \eqref{alpha-interp}.

З формули \eqref{R-interp-F-gen} випливає такий результат для узагальнених просторів Соболєва:

\smallskip

\textbf{Наслідок 2.} \it
Додатково до умови теореми~$1$ припустимо, що функціональні параметри $\alpha_{0}$ і $\alpha_{1}$ належать до $C^{\infty}([1,\infty))$ і задовольняють умову \eqref{cond-for-H-gen}, де замість $\alpha$ узято $\alpha_{0}$ або $\alpha_{1}$. Крім того, припустимо, що $1<p<\infty$. Тоді
\begin{equation}\label{R-interp-H-gen}
(H^{\alpha_0}_{p}(\mathbb{R}^{n}),
H^{\alpha_1}_{p}(\mathbb{R}^{n}))_{\gamma,q}=
B^{\alpha}_{p,q}(\mathbb{R}^{n})
\end{equation}
з точністю до еквівалентності норм. \rm

\smallskip

Формули \eqref{R-interp-B-gen} і \eqref{R-interp-H-gen} доведено в \cite[cс. 193--194, теореми 12 і 13]{Merucci84}. Формулу \eqref{R-interp-F-gen} доведено в \cite[c.~166, теорема~5.3, формула~(2)]{CobosFernandez88} для $p,q_0,q_1\in(1,\infty)$ (там же отримано і формулу \eqref{R-interp-B-gen}). У випадку, коли $p=1$, формула \eqref{R-interp-F-gen} доводиться аналогічно.

З теореми~4 випливає, що кожний простір  $B^{\alpha}_{p,q}(\mathbb{R}^{n})$ отримується дійсною інтерполяцією з функціональним параметром класичних просторів Нікольського\,--\,Бєсова. А~саме, правильне таке твердження:

\smallskip

\textbf{Теорема 5.} \it Нехай $\alpha\in\mathrm{RO}$ і $p,q\in[1,\infty]$. Довільно виберемо дійсні числа $s_0<\sigma_0(\alpha)$ і $s_1>\sigma_1(\alpha)$. Візьмемо неперервну функцію $\widetilde{\alpha}\in\mathrm{RO}$ таку, що $\widetilde{\alpha}(1)=1$ і $\widetilde{\alpha}(t)\asymp\alpha(t)$ при $t\geq1$ (це завжди можна зробити). Покладемо
\begin{equation}\label{def-gamma}
\gamma(t):=\left\{
\begin{array}{ll}
t^{-s_0/(s_1-s_0)}\,\widetilde{\alpha}(t^{1/(s_1-s_0)}),
&\hbox{якщо}\;\;t\geq1;\\
\frac{1}{\gamma(1/t)},&\hbox{якщо}\;\;0<t<1.
\end{array}
\right.
\end{equation}
Тоді функція $\gamma$ належить до класу $\mathfrak{B}$, задовольняє умову \eqref{Boyd-indexes} і з точністю до еквівалентності норм
\begin{equation}\label{B-gen-R-interp}
B^{\alpha}_{p,q}(\mathbb{R}^{n})=
(B^{s_0}_{p,q_{0}}(\mathbb{R}^{n}),
B^{s_1}_{p,q_{1}}(\mathbb{R}^{n}))_{\gamma,q}
\end{equation}
які б не були $q_{0},q_{1}\in[1,\infty]$. \rm

\smallskip

Цей результат наведено в статті \cite[с.~205, твердження~7]{Almeida05} (там використовується множина $\mathfrak{B}$ замість $\mathrm{RO}$ та припускається, що $s_0>\sigma_0(\alpha)$ і $s_1<\sigma_1(\alpha)$). Його окремий випадок розглянуто в \cite[с.~47, формула (2.10)]{CaetanoMoura04} з посиланням на роботи \cite[с.~194]{Merucci84} і \cite[c.~166]{CobosFernandez88}, цитовані вище стосовно теореми~4.

Таким чином, теорема 5 дозволяє отримати усі простори $B^{\alpha}_{p,q}(\mathbb{R}^{n})$, де $\alpha\in\mathrm{RO}$ і $p,q\in[1,\infty]$, за допомогою інтерполяції класичних просторів Нікольського\,--\,Бєсова. Важливість цього факту і корисні наслідки з нього обговорені для (істотно) більш вузького класу параметрів $\alpha$ у монографії Ґ.~Трібеля \cite[с.~58, п.~8]{Triebel10}.

Стосовно формули \eqref{def-gamma} зауважимо, що функція $\gamma_1:=\gamma\!\upharpoonright\mathbb{R}_{+}\!$ належить до класу $\mathrm{RO}$, а її індекси Матушевської
\begin{equation}
\sigma_0(\gamma_1)=\frac{-s_0}{s_1-s_0}+\frac{\sigma_0(\alpha)}{s_1-s_0}>0
\quad\mbox{і}\quad
\sigma_1(\gamma_1)=\frac{-s_0}{s_1-s_0}+\frac{\sigma_1(\alpha)}{s_1-s_0}<1
\end{equation}
за умовою теореми~5. Отже, згідно з \eqref{Matuszewska-Boyd} маємо
\begin{equation}
\nu_0(\overline{\gamma})=\nu_0(\overline{\gamma_1^{\star}})=
\sigma_0(\gamma_1)>0\quad\mbox{і}\quad
\nu_1(\overline{\gamma})=\nu_1(\overline{\gamma_1^{\star}})=
\sigma_1(\gamma_1)<1,
\end{equation}
тобто функція $\gamma$ задовольняє умову \eqref{Boyd-indexes} й тому дійсний метод $(\cdot,\cdot)_{\gamma,q}$ є інтерполяційним функтором.
Формула \eqref{B-gen-R-interp} є окремим випадком формули \eqref{R-interp-B-gen}, де $\alpha_0(t)\equiv t^{s_0}$ і $\alpha_1(t)\equiv t^{s_1}$ та узято $\widetilde{\alpha}$ замість $\alpha$; при цьому рівність \eqref{alpha-R-interp-gamma} випливає з \eqref{def-gamma}.

Окремим і важливим випадком інтерполяційної формули \eqref{B-gen-R-interp} є рівність
\begin{equation}\label{Fp-gen-R-interp}
F^{\alpha}_{p,p}(\mathbb{R}^{n})=
(F^{s_0}_{p,p}(\mathbb{R}^{n}),
F^{s_1}_{p,p}(\mathbb{R}^{n}))_{\gamma,p},
\end{equation}
де $1\leq p<\infty$ (згадаємо \eqref{B=F-gen}). З огляду на це є корисним такий наслідок з теореми~2.

\smallskip

\textbf{Теорема 6.} \it Нехай $\alpha\in\mathrm{RO}$ і $p,q\in[1,\infty)$, причому $p\neq q$. Виберемо число $r$ таке, що $p<r<q$ або $q<r<p$, і означимо число $\theta\in(0,1)$ за формулою
\begin{equation}\label{qpr-interp}
\frac{1}{q}=\frac{1-\theta}{p}+\frac{\theta}{r}.
\end{equation}
Крім того, довільно виберемо число $s_1\in\mathbb{R}$ і означимо функцію $\alpha_0\in\mathrm{RO}$ за формулою~\eqref{alpha-0-interp}. Тоді
\begin{equation}\label{F-gen-C-interp-F}
F^{\alpha}_{p,q}(\mathbb{R}^{n})=[F^{\alpha_0}_{p,p}(\mathbb{R}^{n}),
F^{s_1}_{p,r}(\mathbb{R}^{n})]_{\theta}
\end{equation}
з точністю до еквівалентності норм. \rm

\smallskip

Таким чином, теорема 6 разом з теоремою 5 (розглянутою для $q=q_0=q_1=p$) дозволяє отримати усі простори $F^{\alpha}_{p,q}(\mathbb{R}^{n})$, де $\alpha\in\mathrm{RO}$ і $p,q\in[1,\infty)$, за допомогою повторної інтерполяції класичних просторів Лізоркіна\,--\,Трібеля. Цей результат має перевагу над інтерполяційною формулою \eqref{F-gen-C-interp}, оскільки у ньому є незмінним основний показник сумовності~$p$. Остання обставина буває корисною, коли доводиться інтерполювати лінійні оператори, обмеженість яких на класичних просторах Лізоркіна\,--\,Трібеля залежить від $p$ (наприклад, оператори слідів на краю многовиду та оператори еліптичних крайових задач).

Зауважимо, що, мабуть, уперше дійсна інтерполяція узагальнених (взагалі кажучи, анізотропних) $B$- і $F$-просторів була розглянута Ґ.~Трібелем у статті \cite[с.~239, теорема 4.2/2]{Triebel77III} у випадку числових параметрів інтерполяції. Втім, як зазначалося вище, класи просторів, уведених Ґ.~Трібелем у цій статті, не містять просторів, розглянутих нами, за виключенням узагальнених $H$-просторів. Пізніше з'явилася серія робіт А.~Г.~Багдасаряна \cite{Bagdasaryan92, Bagdasaryan96, Bagdasaryan97, Bagdasaryan01, Bagdasaryan04, Bagdasaryan05, Bagdasaryan10}, у якій систематично застосовується і досліджується дійсна інтерполяція (також з числовими параметрами) анізотропних узагальнених $H$-, $B$- і $F$-просторів. Зокрема, у його статтях \cite{Bagdasaryan04, Bagdasaryan05, Bagdasaryan10} узагальнені $B$-простори уведені за допомогою дійсної інтерполяції узагальнених $H$-просторів. Базові для нас роботи Ч.~Meруччі \cite{Merucci84} і Ф.~Кобоса, Д.~Л.~Фернандеза \cite{CobosFernandez88} є піонерським у застосуванні дійсної інтерполяції з функціональним параметром до просторів узагальненої гладкості. Їх результати використовували і узагальнювали у низці робіт (їх кілька десятків) вже у цьому столітті; серед них відмітимо статті \cite{CaetanoMoura04, Almeida05, Almeida09, AlmeidaCaetano11, CaetanoLopes11, CobosDominguez14, CobosDominguez15, BesoyCobos18}. Деяка частина цих робіт присвячена логарифмічним просторам узагальненої гладкості, для яких показник регулярності має (у наших позначеннях) вигляд $\alpha(t)\equiv t^{s}(1+\log t)^{r}$, де $s,r\in\mathbb{R}$; див., наприклад, монографію \cite[пп. 1.3, 3.4]{Triebel10}, статті  \cite{CobosDominguez14, CobosDominguez15}, великий за обсягом препринт \cite{DominguezTikhonov18} та наведену там бібліографію. Мабуть, ці простори утворюють найбільш важливий для застосувань клас просторів узагальненої гладкості.
Для дослідження логарифмічних $B$-просторів використовують дійсну інтерполяцію з (першим) функціональним показником вигляду $\gamma(t)=t^{\theta}(1+\log|t|)^{r}$, якщо $0<t<1$, і $\gamma(t)=t^{\theta}(1+\log|t|)^{l}$, якщо $t\geq1$, де $0<\theta<1$ та $r,l\in\mathbb{R}$ (розглядаються також і граничні випадки $\theta=0$ та $\theta=1$). Такій інтерполяції присвячені, наприклад, статті \cite{Doktorskii91, EvansOpic00, EvansOpic02, BesoyCobos18}. У новітній роботі \cite[с.~19, наслідок~8.6]{LoosveldtNicolay19} показано, що простори $B^{\alpha}_{p,q}(\mathbb{R}^{n})$, де $\alpha\in\mathrm{RO}$, отримуються деяким методом дійсної інтерполяції (з функціональним параметром) класичних просторів Соболєва, причому охоплено і граничні випадки $p=1$ і $p=\infty$. Цей метод уведено у цій же роботі й він, взагалі кажучи, відрізняється від традиційного методу, використаного в \cite{Merucci84, CobosFernandez88, CaetanoMoura04, Almeida05}.

Порівнюючи комплексний і дійсний методи інтерполяції, природно поставити таке питання: чи можливо у який-небудь спосіб узагальнити комплексний метод на випадок функціонального параметра інтерполяції так, щоб отримати простори узагальненої гладкості шляхом інтерполяції цим методом класичних просторів. Тут особливо цікавим є випадок $H$- і $F$-просторів (бо для $B$-просторів достатньо вже відомого методу дійсної інтерполяції з функціональним параметром). Наскільки нам відомо, поки не знайдено таких узагальнень, які дають усі простори $F^{\alpha}_{p,q}(\mathbb{R}^{n})$, де $\alpha\in\mathrm{RO}$, шляхом інтерполяції класичних просторів Лізоркіна\,--\,Трібеля, навіть у випадку $q=2$ (тобто узагальнених просторів Соболєва). Є лише декілька узагальнень методу комплексної інтерполяції \cite{Lions60/61, Schechter67, CarroCerda90, FanKaijser94, Fan94}, які можна трактувати як інтерполяцію з функціональним параметром; проте,  використані у цих роботах класи параметрів інтерполяції занадто вузькі для наших цілей. Серед цих робіт дві \cite{Schechter67, CarroCerda90} дають застосування уведених методів комплексної інтерполяції до просторів узагальненої гладкості, а саме, $H$-просторів.

\medskip

\textbf{4. Простори узагальненої гладкості на многовидах.} Узагальнені $B$- і $F$-простори, розглянуті вище, є інваріантними відносно $C^{\infty}$-дифеоморфізмів простору $\mathbb{R}^{n}$ на себе, які є тотожними відображеннями в околі нескінченності (такі дифеоморфізми називаємо локальними). Це випливає на підставі теорем 5 і 6 з інваріантності класичних $B$- і $F$-просторів відносно цих дифеоморфізмів (див. \cite[с.~247, п.~2.10.2, теорема]{Triebel86}). Ця властивість дозволяє коректно означити $B$- і $F$-простори на компактних $C^{\infty}$-многовидах за допомогою скінченних наборів локальних карт, які покривають многовид. Тут коректність розуміється у тому сенсі, що уведені так простори не залежать (з точністю до еквівалентності норм) від вибору набору локальних карт. Нагадаємо, що $C^{\infty}$-дифеоморфізмом простору $\mathbb{R}^{n}$ називають будь-яке взаємно однозначне відображення $\pi:\mathbb{R}^{n}\leftrightarrow\mathbb{R}^{n}$, яке нескінченно диференційовне на $\mathbb{R}^{n}$ разом з оберненим відображенням $\pi^{-1}$. Локальний $C^{\infty}$-дифеоморфізм $\pi$ породжує лінійний оператор заміни змінної $\pi:f\mapsto f\circ\pi$, де $f\in\mathcal{S}'(\mathbb{R}^{n})$ (якщо $f$~--- класична функція, то $(f\circ\pi)(x)=f(\pi(x))$). Вказана інваріантність означає, що цей оператор є ізоморфізмом (тобто взаємно однозначним і взаємно неперервним лінійним відображенням) $B$- і $F$-просторів на себе. Відмітимо, що вона випливає з \cite[с.~175, теорема~3]{Kalyabin85} для просторів $B^{\alpha}_{p,q}(\mathbb{R}^{n})$ і $F^{\alpha}_{p,q}(\mathbb{R}^{n})$ у випадку, коли функція $\alpha\in\mathrm{RO}$ задовольняє умову $\sigma_0(\alpha)>0$.

Питання про інваріантність деяких узагальнених $B$- і $F$-просторів відносно дифеоморфізмів (тобто про локалізацію цих просторів) і, як наслідок, коректність їх означення на нескінченно гладких многовидах поставлено у монографії Ґ.~Трібеля \cite[с.~57, пп.~4 і~6]{Triebel10}.

Для того, щоб означити узагальнені $B$- і $F$-простори на многовиді з краєм, потрібні їх аналоги на евклідовому півпросторі $\mathbb{R}^{n}_{+}$, оскільки, не обмежуючи загальності, можна вважати, що локальні карти в околах крайових точок многовиду здійснюють гомеоморфізми між цими околами і $\mathbb{R}^{n}_{+}$. Такі аналоги означаються у стандартний спосіб як простори звужень на $\mathbb{R}^{n}_{+}$ розподілів з відповідних просторів на $\mathbb{R}^{n}$.  Цей підхід застосовують і для просторів, заданих на довільній відкритій непорожній множині $G\subset\mathbb{R}^{n}$ (див., наприклад, монографії \cite[c.~384, п.~4.2.1, означення~1]{Triebel80} і \cite[c.~272, п.~3.2.2, означення~1]{Triebel86}, де розглянуті класичні $B$- і $F$-простори). Наведемо відповідне означення для довільного банахового простору $\Upsilon(\mathbb{R}^{n})$, який є лінійним многовидом у просторі $\mathcal{S}'(\mathbb{R}^{n})$ і неперервно вкладається у цей простір.

\smallskip

\textbf{Означення 5.} Лінійний простір $\Upsilon(G)$ складається із звужень в $G$ усіх розподілів $f\in\Upsilon(\mathbb{R}^{n})$ і наділений нормою
\begin{equation}
\|u\|_{\Upsilon(G)}:=\inf\bigl\{\|f\|_{\Upsilon(\mathbb{R}^{n})}:
f\in\Upsilon(\mathbb{R}^{n}),\;\;f=u\;\;\mbox{в}\;\;G\bigr\},
\end{equation}
де $u\in\Upsilon(G)$.

\smallskip

Простір $\Upsilon(G)$ повний, тобто банахів, оскільки він є факторпростором повного простору $\Upsilon(\mathbb{R}^{n})$ за його підпростором
\begin{equation}
\bigl\{f\in\Upsilon(\mathbb{R}^{n}):
\mathrm{supp}\,f\subset\mathbb{R}^{n}\setminus G\bigr\}.
\end{equation}
Замкненість останнього випливає з неперервності вкладення $\Upsilon(\mathbb{R}^{n})$ в $\mathcal{S}'(\mathbb{R}^{n})$. З наведеного означення негайно випливає, що лінійний оператор звуження $R_{G}:f\mapsto f\!\upharpoonright\!G$ розподілу $f\in\mathcal{S}'(\mathbb{R}^{n})$ на множину $G$ є сюр'єктивним обмеженим оператором $R_{G}:\Upsilon(\mathbb{R}^{n})\to\Upsilon(G)$.

Узявши тут $\Upsilon(\mathbb{R}^{n}):=B^{\alpha}_{p,q}(\mathbb{R}^{n})$ або $\Upsilon(\mathbb{R}^{n}):=F^{\alpha}_{p,q}(\mathbb{R}^{n})$, отримуємо відповідно банахові простори $B^{\alpha}_{p,q}(G)$ і $F^{\alpha}_{p,q}(G)$, де $\alpha\in\mathrm{RO}$, а $p,q\in[1,\infty]$ для $B$-простору, або $p,q\in[1,\infty)$ для $F$-простору.

\smallskip

\textbf{Теорема 7.} \it Припустимо, що $G=\mathbb{R}^{n}_{+}$ або $G$~--- обмежена область в $\mathbb{R}^{n}$ з ліпшіцевою межею. Тоді усі твердження теорем $2$--$6$ зберігають силу, якщо у них замінити $\mathbb{R}^{n}$ на~$G$. \rm

\smallskip

Потрібні інтерполяційні формули для просторів, заданих на $G$, виводяться з їх аналогів для просторів, заданих на $\mathbb{R}^{n}$, за допомогою обмеженого лінійного оператора $T:\Upsilon(G)\to\Upsilon(\mathbb{R}^{n})$ такого, що  $R_{G}T$~--- тотожний оператор на $\Upsilon(G)$. Цей оператор продовжує розподіл $u\in\Upsilon(G)$ до розподілу, заданого на $\mathbb{R}^{n}$. Тут $\Upsilon(\mathbb{R}^{n})$ позначає простори, що утворюють інтерполяційну пару, а оператор $T$ діє в один і той же спосіб на перетині цих просторів. Тоді простори $\Upsilon(G)$ є ретрактами просторів $\Upsilon(\mathbb{R}^{n})$ (з відповідними ретракцією $R_{G}$ і ретракцією $T$) й тому інтерполяція просторів $\Upsilon(G)$ зводиться до інтерполяції просторів $\Upsilon(\mathbb{R}^{n})$ (порівняти, наприклад, з \cite[c. 151--152, доведення теореми 3.2]{MikhailetsMurach10}). Вказаний оператор $T$ побудовано для усіх класичних просторів $B^{s}_{p,q}(G)$ і $F^{s}_{p,q}(G)$ В.~C.~Ричковим \cite[с.~253, теорема~4.1]{Rychkov99} (його називають універсальним, див., наприклад, монографію \cite[с.~66]{Triebel06}). У випадку $G=\mathbb{R}^{n}_{+}$ і $|s|\leq n$, де $n\in\mathbb{N}$, оператор $T$, залежний від $n$, розглянуто, наприклад, в \cite[п.~4.5.2, с.~225, наслідок]{Triebel92}. Спочатку доводиться аналог теореми~5 для просторів на $G$, оскільки у ній інтерполюється пара класичних просторів. З неї випливає, що $T$ є обмеженим оператором з $B^{\alpha}_{p,q}(G)$ в $B^{\alpha}_{p,q}(\mathbb{R}^{n})$ для довільного $\alpha\in\mathrm{RO}$. Це дозволяє довести аналог теореми~6, з якої випливає, що $T$ є також обмеженим оператором з $F^{\alpha}_{p,q}(G)$ в $F^{\alpha}_{p,q}(\mathbb{R}^{n})$. Тепер можна довести аналоги теорем~2 і~4, а теорема~3 є безпосереднім наслідком теореми~2.

\smallskip

\textbf{Наслідок 3.} \it За умови теореми $7$ універсальний лінійний оператор продовження Ричкова є обмеженим оператором
$B^{\alpha}_{p,q}(G)\to B^{\alpha}_{p,q}(\mathbb{R}^{n})$ для усіх $\alpha\in\mathrm{RO}$ і $p,q\in[1,\infty]$, а також є обмеженим оператором
$F^{\alpha}_{p,q}(G)\to F^{\alpha}_{p,q}(\mathbb{R}^{n})$ для усіх $\alpha\in\mathrm{RO}$ і $p,q\in[1,\infty)$. \rm

\smallskip

Це випливає з теореми 7. У монографії Ґ.~Трібеля \cite[с.~58, п.~8]{Triebel10} теорему~7 разом із наслідком~3 сформульовано для логарифмічних узагальнених $B$-просторів.

Перейдемо до просторів узагальненої гладкості на $C^{\infty}$-многовидах. Нехай $M$~--- нескінченно гладкий орієнтовний компактний многовид вимірності $n\geq1$ з краєм $\partial M$. Нехай
$M^{\circ}:=M\setminus\partial M$. Допускається також випадок, коли  $\partial M=\varnothing$, тобто коли $M=M^{\circ}$~--- замкнений многовид. Як і в \cite[c.~636]{Hermander87}, позначимо через $\overline{\mathcal{D}}\,'(M^{\circ})$ лінійний топологічний простір усіх продовжуваних розподілів на $M^{\circ}$. (Якщо многовид $M$ замкнений, то $\overline{\mathcal{D}}\,'(M^{\circ})$~--- це топологічний простір  $\mathcal{D}'(M)$ усіх розподілів на $M$.)

Виберемо який-небудь скінченний атлас $\pi_{j}:\overline{\Pi}_{j}\leftrightarrow
U_{j}$, де $j=1,\ldots,\varkappa$, з $C^{\infty}$-структури на $M$. Тут усі $U_{j}$~--- відкриті (у топології простору $M$) множини, які утворюють скінченне покриття многовиду $M$, а кожне $\Pi_{j}$ позначає або $\mathbb{R}^{n}$ (якщо $U_{j}\cap\partial M=\varnothing$) або $\mathbb{R}^{n}_{+}$ (у противному випадку). Як звичайно,  $\overline{\Pi}_{j}$~--- замикання множини $\Pi_{j}$ в $\mathbb{R}^{n}$, тобто $\overline{\Pi}_{j}$~--- або $\mathbb{R}^{n}$, або $\overline{\mathbb{R}^{n}_{+}}:=
\{(x',x_n):x'\in\mathbb{R}^{n-1},x_n\geq0\}$, відповідно. (Якщо многовид $M$ замкнений, то усі $\Pi_{j}=\overline{\Pi}_{j}=\mathbb{R}^{n}$.) Крім того, виберемо яке-небудь розбиття одиниці $\chi_{j}\in
C^{\infty}(M)$, де $j=1,\ldots,\varkappa$, на $M$, яке задовольняє умову $\mathrm{supp}\,\chi_{j}\subset U_{j}$.

За допомогою банахових просторів $\Upsilon(\mathbb{R}^{n})$ і $\Upsilon(\mathbb{R}^{n}_{+})$, розглянутих вище даємо

\smallskip

\textbf{Означення 6.} Лінійний простір $\Upsilon(M^{\circ})$ складається з усіх розподілів $g\in\overline{\mathcal{D}}\,'(M^{\circ})$ таких, що
$(\chi_{j}g)\circ\pi_{j}\in\Upsilon(\Pi_{j})$ для кожного
$j\in\{1,\ldots,\varkappa\}$; тут $(\chi_{j}g)\circ\pi_{j}$~--- зображення розподілу $\chi_{j}g$ у локальній карті $\pi_{j}$. Цей простір наділено нормою
\begin{equation}
\|g\|_{\Upsilon(M^{\circ})}:=\sum_{j=1}^{\varkappa}
\|(\chi_{j}g)\circ\pi_{j}\|_{\Upsilon(\Pi_{j})}.
\end{equation}

Узявши тут усі $\Upsilon(\Pi_{j}):=B^{\alpha}_{p,q}(\Pi_{j})$ або усі $\Upsilon(\Pi_{j}):=F^{\alpha}_{p,q}(\Pi_{j})$, отримуємо відповідно банахові простори $B^{\alpha}_{p,q}(M^{\circ})$ і $F^{\alpha}_{p,q}(M^{\circ})$, де $\alpha\in\mathrm{RO}$, а $p,q\in[1,\infty]$ для $B$-простору, або $p,q\in[1,\infty)$ для $F$-простору.

\smallskip

\textbf{Теорема 8.} \it Уведені простори $B^{\alpha}_{p,q}(M^{\circ})$ і $F^{\alpha}_{p,q}(M^{\circ})$ є повними та з точністю до еквівалентності норм не залежать від зазначеного вибору атласу многовиду і розбиття одиниці. \rm

\smallskip

\textbf{Теорема 9.} \it Усі твердження теорем $2$--$6$ зберігають силу, якщо у них замінити $\mathbb{R}^{n}$ на~$M^{\circ}$. \rm

\smallskip

Обговоримо доведення цих теорем. Згідно з означенням~6, лінійне відображення розпрямлення
\begin{equation}\label{T-distribution}
T:g\mapsto
\bigl((\chi_{1}g)\circ\pi_{1},\ldots,
(\chi_{\varkappa}g)\circ\pi_{\varkappa}\bigr),
\quad\mbox{де}\;\;g\in\overline{\mathcal{D}}\,'(M^{\circ}),
\end{equation}
є ізометричним оператором
\begin{equation}\label{T-Upsilon}
T:\Upsilon(M^{\circ})\to\prod_{j=1}^{\varkappa}\Upsilon(\Pi_{j}).
\end{equation}
Для $T$ існує лінійне відображення склеювання
\begin{equation}\label{K-distribution}
K:\prod_{j=1}^{\varkappa}\mathcal{S}'(\Pi_{j})
\to\overline{\mathcal{D}}\,'(M^{\circ})
\end{equation}
таке, що $KTg=g$ для довільного $g\in\overline{\mathcal{D}}\,'(M^{\circ})$. Тут, звичайно, $\mathcal{S}'(\mathbb{R}^{n}_{+})$~--- лінійний топологічний простір звужень на $\mathbb{R}^{n}_{+}$ усіх розподілів з простору $\mathcal{S}'(\mathbb{R}^{n})$ (це стосується випадку $\Pi_{j}=\mathbb{R}^{n}_{+}$). Відображення $K$ неважко побудувати; див., наприклад, \cite[c.~85]{MikhailetsMurach10} у випадку замкненого многовиду. Якщо простір $\Upsilon(\mathbb{R}^{n})$ інваріантний відносно локальних $C^{\infty}$-дифеоморфізмів простору $\mathbb{R}^{n}$ на себе, то відображення склеювання є обмеженим оператором
\begin{equation}\label{K-Upsilon}
K:\prod_{j=1}^{\varkappa}\Upsilon(\Pi_{j})\to\Upsilon(M^{\circ}).
\end{equation}
Як зазначалося вище, розглянуті простори $B^{\alpha}_{p,q}(\mathbb{R}^{n})$ і $F^{\alpha}_{p,q}(\mathbb{R}^{n})$ інваріантні відносно цих дифеоморфізмів. Тому простори $B^{\alpha}_{p,q}(M^{\circ})$ і $F^{\alpha}_{p,q}(M^{\circ})$ є ретрактами відповідно просторів
\begin{equation}
\prod_{j=1}^{\varkappa}B^{\alpha}_{p,q}(\Pi_{j})\quad\mbox{і}\quad
\prod_{j=1}^{\varkappa}F^{\alpha}_{p,q}(\Pi_{j})
\end{equation}
(з ретракцією $K$ і коретракцією $T$). Отже, оскільки останні два простори повні, то і перші два також повні. Крім того, з теорем 2--6 (для $\Pi_{j}=\mathbb{R}^{n}$) і теореми~7 (для $\Pi_{j}=\mathbb{R}^{n}_{+}$) випливають їх аналоги для просторів на $M^{\circ}$ (див. \cite[c.~84--87, доведення теореми~2.2]{MikhailetsMurach10}). При цьому використовується така загальна властивість інтерполяційних методів (функторів): інтерполяція прямого добутку пар банахових просторів збігається (з точністю до еквівалентності норм) з прямим добутком просторів, отриманих інтерполяцією цих пар. Незалежність (з точністю до еквівалентності норм) просторів $B^{\alpha}_{p,q}(M^{\circ})$ і $F^{\alpha}_{p,q}(M^{\circ})$  від зазначеного вибору атласу многовиду і розбиття одиниці випливає з  незалежності класичних $B$- і $F$-просторів, заданих на $M^{\circ}$, від цього вибору та теореми~9 (аналогів формул \eqref{B-gen-R-interp} і \eqref{F-gen-C-interp-F} для просторів на $M^{\circ}$). Остання незалежність доведена в \cite[с.~272, п.~3.2.3, твердження (ii), (iii)]{Triebel86} у випадку $\partial M=\varnothing$. У противному випадку доведення аналогічне. Втім, неважко безпосередньо довести вказану незалежність просторів $B^{\alpha}_{p,q}(M^{\circ})$ і $F^{\alpha}_{p,q}(M^{\circ})$, скориставшись їх інваріантністю відносно локальних $C^{\infty}$-дифеоморфізмів.

Якщо $G$~--- обмежена область в $\mathbb{R}^{n}$ з межею $\partial G\in C^{\infty}$, то її замикання $\overline{G}$ є прикладом зазначеного многовиду~$M$ (тоді $G=M^{\circ}$). У цьому випадку маємо різні означення 5 і 6 кожного з просторів $B^{\alpha}_{p,q}(G)$ і $F^{\alpha}_{p,q}(G)$.

\smallskip

\textbf{Теорема 10.} \it У вказаному випадку означення $5$ і $6$ еквівалентні у тому сенсі, що вони уводять одні і ті ж самі узагальнені $B$- і $F$-простори з точністю до еквівалентності норм. \rm

\smallskip

З інваріантності просторів $B^{\alpha}_{p,q}(\mathbb{R}^{n})$ і $F^{\alpha}_{p,q}(\mathbb{R}^{n})$ відносно локальних $C^{\infty}$-дифеомор\-фіз\-мів простору $\mathbb{R}^{n}$ на себе випливає, що відображення \eqref{T-distribution} і \eqref{K-distribution} задають відповідно обмежені оператори \eqref{T-Upsilon} і \eqref{K-Upsilon}, де $\Upsilon(M^{\circ})$ позначає один з просторів $B^{\alpha}_{p,q}(G)$ і $F^{\alpha}_{p,q}(G)$, уведений в означенні~5, а $\Upsilon(\Pi_{j})$ позначає їх аналоги на $\Pi_{j}$. Крім того, як вказано вище, відображення \eqref{T-distribution} і \eqref{K-distribution} задають відповідно обмежені оператори \eqref{T-Upsilon} і \eqref{K-Upsilon}, де $\Upsilon(M^{\circ})$ позначає один з просторів $B^{\alpha}_{p,q}(G)$ і $F^{\alpha}_{p,q}(G)$, уведений в означенні~6, а $\Upsilon(\Pi_{j})$ є тим же самим простором. Звідси негайно випливає висновок теореми~10. Інший підхід до доведення цієї теореми: скористатися теоремою~9 стосовно інтерполяційних формул  \eqref{B-gen-R-interp} і \eqref{F-gen-C-interp-F} та тим фактом, що висновок теореми~10 правильний для класичних $B$- і $F$-просторів.

Теорема 10 дає позитивну відповідь на питання Ґ.~Трібеля \cite[с.~57, п.~4]{Triebel10} про можливість локалізації деяких узагальнених $B$- і $F$-просторів на $G$, зокрема, логарифмічних просторів.

\medskip

\textbf{5. Застосування до еліптичних операторів.} Розглянемо застосування узагальнених $B$- і $F$-просторів на многовидах до еліптичних псевдодиференціальних операторів (ПДО). Припустимо, що нескінченно гладкий компактний многовид $M$ є замкненим, тобто його край $\partial M=\varnothing$. Позначимо через $\Psi^{r}_{\mathrm{ph}}(M)$ множину усіх класичних (тобто поліоднорідних) ПДО порядку $r\in\mathbb{R}$, заданих на $M$ (див. відповідні означення в  монографії \cite[пп. 1.4.1 і 2.2.1]{MikhailetsMurach10}  або огляді \cite[\S\S~1 і~2]{Agranovich90}). Важливим прикладом ПДО класу $\Psi^{r}_{\mathrm{ph}}(M)$, де $r\in\mathbb{N}$, є довільний лінійний диференціальний оператор на $M$ порядку $r$ з нескінченного гладкими коефіцієнтами на $M$.

Кожний ПДО $A\in\Psi^{r}_{\mathrm{ph}}(M)$ є неперервним лінійним оператором на топологічному просторі $\mathcal{D}'(M)$ усіх розподілів на~$M$. Будемо розглядати звуження цього ПДО на узагальнені $B$- і $F$-простори на~$M$; ці звуження позначаємо також через~$A$.

\smallskip

\textbf{Лема 1.} \it Нехай $r\in\mathbb{R}$ і $A\in\Psi^{r}_{\mathrm{ph}}(M)$. Тоді ПДО $A$ є обмеженим оператором на парі просторів
\begin{equation}\label{PsDO-B}
A:B^{\alpha}_{p,q}(M)\to B^{\alpha\varrho^{-r}}_{p,q}(M)
\end{equation}
для усіх $\alpha\in\mathrm{RO}$ і $p,q\in[1,\infty]$, а також на парі просторів
\begin{equation}\label{PsDO-F}
A:F^{\alpha}_{p,q}(M)\to F^{\alpha\varrho^{-r}}_{p,q}(M)
\end{equation}
для усіх $\alpha\in\mathrm{RO}$ і $p,q\in[1,\infty)$. \rm

\smallskip

Тут і надалі у показнику регулярності просторів використовуємо функціональний параметр $\varrho(t):=t$ аргументу $t\geq1$. Отже, $\alpha\varrho^{-r}$ у лемі~1 позначає функцію $\alpha(t)t^{-r}$ цього аргументу.

Якщо функція $\alpha$ степенева, ця лема випливає з \cite[с.~258, теорема п.~6.2.2]{Triebel92}. Вказана теорема стосується ПДО порядку $r=0$, заданих на $\mathbb{R}^{n}$. Випадок довільного $r$ зводиться до випадку $r=0$ множенням цих ПДО на ПДО
\begin{equation}
J^{-r}:w\mapsto
\mathcal{F}^{-1}[(\langle\xi\rangle)^{-r}(\mathcal{F}w)(\xi)],
\quad\mbox{де}\;\;w\in\mathcal{S}'(\mathbb{R}^{n}),
\end{equation}
порядку $-r$. Для ПДО $A$ на $M$ потрібний результат отримаємо, скориставшись його зображенням \cite[п.~2.1, с.~22, формула (2.1.3)]{Agranovich90} за допомогою ПДО на $\mathbb{R}^{n}$. Для довільного $\alpha\in\mathrm{RO}$ лема~1 виводиться з випадку степеневої функції $\alpha$ за допомогою інтерполяції на підставі теореми~9 (використовуються аналоги формул \eqref{B-gen-R-interp} і \eqref{F-gen-C-interp-F} для просторів на~$M$).

Топологічний простір $\mathcal{D}'(M)$ усіх розподілів на $M$ є спряженим до простору $C^{\infty}(M)$ відносно скалярного добутку у гільбертовому просторі $L_{2}(M,dx)$ функцій, квадратично інтегровних на многовиді $M$ відносно деякої $C^{\infty}$-щільності $dx$, заданої на $M$. Позначимо через $A^{+}$ ПДО, формально спряжений до ПДО $A\in\Psi^{r}_{\mathrm{ph}}(M)$ відносно цього скалярного добутку. Отже, $(A^{+}f,u)_{M}=(f,Au)_{M}$ для будь-яких $f\in\mathcal{D}'(M)$ і $u\in C^{\infty}(M)$, причому $A^{+}\in\Psi^{r}_{\mathrm{ph}}(M)$; тут $(\cdot,\cdot)_{M}$ позначає продовження за неперервністю вказаного скалярного добутку. Покладемо
\begin{gather}\label{kernel-A}
\mathcal{N}:=\{u\in C^{\infty}(M): Au=0\;\mbox{на}\;M\},\\
\mathcal{N}^{+}:=\{v\in C^{\infty}(M): A^{+}v=0\;\mbox{на}\;M\}.
\label{cokernel-A}
\end{gather}

\textbf{Теорема 11.} \it Нехай $r\in\mathbb{R}$ і $A\in\Psi^{r}_{\mathrm{ph}}(M)$. Припустимо, що ПДО $A$ є еліптичним на многовиді~$M$. Тоді обмежені лінійні оператори \eqref{PsDO-B} і \eqref{PsDO-F} є нетеровими для вказаних значень параметрів $\alpha$, $p$ і $q$. Ядро кожного з цих операторів збігається з $\mathcal{N}$, а область значень складається з усіх розподілів $f\in B^{\alpha\varrho^{-r}}_{p,q}(M)$ для оператора \eqref{PsDO-B} або $f\in F^{\alpha\varrho^{-r}}_{p,q}(M)$ для оператора \eqref{PsDO-F} таких, що $(f,v)_{M}=0$ для довільного $v\in\mathcal{N}^{+}$. Індекс кожного з цих операторів дорівнює $\dim\mathcal{N}-\dim\mathcal{N}^{+}$ (та є нулем, якщо $\dim M\geq2$). \rm

\smallskip

Ця теорема випливає з випадку степеневої функції $\alpha$ за допомогою інтерполяції на підставі теореми~9 і теореми про інтерполяцію нетерових операторів (остання аналогічна \cite[с.~35, теорема~1.7]{MikhailetsMurach10} і справедлива для довільних інтерполяційних функторів). Крім того, ця теорема нескладно виводиться з того факту \cite[с.~26, теорема 2.2.4]{Agranovich90}, що для еліптичного ПДО $A$ існує параметрикс, тобто такий еліптичний ПДО $R\in\Psi^{-r}_{\mathrm{ph}}(M)$, що
\begin{equation}\label{parametrix-equa}
RA=I+T_1\quad\mbox{і}\quad AR=I+T_2,
\end{equation}
де $I$~--- тотожний оператор на $\mathcal{D}'(M)$, а $T_1$ і $T_2$~--- деякі ПДО на $M$ порядку $-\infty$, тобто  $T_1,T_2\in\Psi^{l}_{\mathrm{ph}}(M)$ для кожного $l\in\mathbb{R}$. Обговоримо доведення теореми~11 за допомогою формули \eqref{parametrix-equa} у випадку оператора \eqref{PsDO-F}. За лемою~1, маємо обмежений оператор
\begin{equation}\label{parametrix-F}
R:F^{\alpha\varrho^{-r}}_{p,q}(M)\to F^{\alpha}_{p,q}(M)
\end{equation}
та компактні оператори
\begin{gather}\label{T1-comp}
T_1:F^{\alpha}_{p,q}(M)\to F^{\alpha}_{p,q}(M),\\
T_2:F^{\alpha\varrho^{-r}}_{p,q}(M)\to F^{\alpha\varrho^{-r}}_{p,q}(M).
\end{gather}
Компактність, наприклад, оператора \eqref{T1-comp} випливає з обмеженості оператора
\begin{equation}\label{T1-F-F}
T_1:F^{\alpha}_{p,q}(M)\to F^{\alpha\varrho^{k}}_{p,q}(M)
\end{equation}
при $k\geq1$ і вкладень
\begin{equation}\label{embed-FHHF}
F^{\alpha\varrho^{k}}_{p,q}(M)\hookrightarrow H^{s}_{p}(M)\hookrightarrow
H^{s-1}_{p}(M)\hookrightarrow F^{\alpha}_{p,q}(M),
\end{equation}
де число $s$ пов'язане з індексами Матушевської функції $\alpha$ умовами $\sigma_1(\alpha)<s-1$ і $s<k+\sigma_0(\alpha)$; ці вкладення неперервні, причому середнє з них компактне. Отже, оператор \eqref{parametrix-F} є двобічним регуляризатором оператора \eqref{PsDO-F}, що еквівалентно нетеровості останнього. Далі, якщо розподіл $u\in F^{\alpha}_{p,q}(M)$ задовольняє умову $Au=0$ на $M$, то
\begin{equation}\label{u-in-C-infty}
u=-T_1u\in\bigcap_{s\in\mathbb{R}}H^{s}_{p}(M)=C^{\infty}(M)
\end{equation}
на підставі першої рівності \eqref{parametrix-equa}, співвідношення \eqref{T1-F-F} і лівого вкладення \eqref{embed-FHHF}. Отже, $\mathcal{N}$~--- ядро оператора \eqref{PsDO-F}. Якщо розподіл $f\in F^{\alpha\varrho^{-r}}_{p,q}(M)$ задовольняє умову $(f,v)_{M}=0$ для усіх $v\in\mathcal{N}^{+}$, то на підставі другої рівності \eqref{parametrix-equa} і аналога вкладення \eqref{u-in-C-infty} робимо висновок, що функція $T_{2}f\in C^{\infty}(M)$ задовольняє умову $(T_{2}f,v)_{M}=0$ для довільного $v\in\mathcal{N}^{+}$. Тому згідно з \cite[с.~31, теорема 2.3.12]{Agranovich90} і \cite[с.~26, теорема 2.2.6]{Agranovich90} існує функція $u\in C^{\infty}(M)$ така, що $Au=T_{2}f$ на $M$. Отже $A(Rf-u)=f$ за другої рівністю \eqref{parametrix-equa}, де $Rf-u\in F^{\alpha}_{p,q}(M)$, тобто $f$ належить області значень оператора \eqref{PsDO-F}. Звісно, кожний розподіл $f$ з цієї області значень задовольняє умову $(f,v)_{M}=0$ для усіх $v\in\mathcal{N}^{+}$. Із отриманих описів ядра і області визначення оператора \eqref{PsDO-F} випливає зазначена формула його індексу. Як відомо \cite[с.~32]{Agranovich90}, він дорівнює нулю, якщо $\dim M\geq2$.

Версія цієї теореми правильна і для матричних еліптичних ПДО вигляду  $\mathbf{A}=(A_{j,k})_{j,k=1}^{\lambda}$, де кожне $A_{j,k}$~--- класичний ПДО на $M$ порядку $\mathrm{ord}\,A_{j,k}\leq r_{j}+l_{k}$. Тут задано ціле число $\lambda\geq2$, набір дійсних чисел $r_{1},\ldots,r_{\lambda}$ і $l_{1},\ldots,l_{\lambda}$ та припускається, що $\mathbf{A}$ є еліптичним на $M$ за Дуглісом--Ніренбергом для цього набору чисел (див., наприклад, \cite[п.~3.2~b]{Agranovich90}).

\smallskip

\textbf{Теорема 11М.} \it Матричний ПДО $\mathbf{A}$ є нетеровим обмеженим оператором на парі просторів
\begin{equation}\label{PsDO-matrix-B}
\mathbf{A}:\prod_{k=1}^{\lambda}B^{\alpha\varrho^{l_k}}_{p,q}(M)\to
    \prod_{j=1}^{\lambda}B^{\alpha\varrho^{-r_j}}_{p,q}(M)
\end{equation}
для усіх $\alpha\in\mathrm{RO}$ і $p,q\in[1,\infty]$, а також на парі просторів
\begin{equation}\label{PsDO-matrix-F}
\mathbf{A}:\prod_{k=1}^{\lambda}F^{\alpha\varrho^{l_k}}_{p,q}(M)\to
    \prod_{j=1}^{\lambda}F^{\alpha\varrho^{-r_j}}_{p,q}(M)
\end{equation}
для усіх $\alpha\in\mathrm{RO}$ і $p,q\in[1,\infty)$. Ядро і коядро кожного з цих операторів лежать у просторі $(C^{\infty}(M))^{\lambda}$ і разом з індексом не залежать від параметрів $\alpha$, $p$ і~$q$. \rm

\smallskip

З теорем 11 і 11М випливають точні достатні умови локальної регулярності розв'язків еліптичних рівнянь в узагальнених $B$- і $F$-просторах та відповідні апріорні оцінки розв'язків (порівняти, наприклад, з \cite[теореми~8 і~9]{Zinchenko17}).

\medskip

\textbf{6. Застосування до еліптичних крайових задач.} Перейдемо до застосувань узагальнених $B$- і $F$-просторів до еліптичних крайових задач (ЕКЗ). Відмітимо, що питання про дослідження ЕКЗ в узагальнених $B$-просторах було поставлено Ґ.~Трібелем в  \cite[с.~59]{Triebel10}.
Припустимо, що $\dim M\geq2$ і $\partial M\neq\varnothing$, та нагадаємо, що $\nobreak{M^\circ:=M\setminus\partial M}$. Зауважимо, що $\partial M$ є замкненим компактним орієнтовним многовидом класу $C^{\infty}$ і вимірності $n-1$. Отже, на $M^\circ$ і $\partial M$ означені узагальнені $B$- і $F$-простори.

Розглянемо на $M^\circ$ ЕКЗ вигляду
\begin{gather}\label{ep-equ}
Lu=f\quad\mbox{на}\;\,M^\circ,\\
B_{j}u=g_{j}\quad\mbox{на}\;\,\partial M, \label{ep-bound-cond}
\quad j=1,...,r
\end{gather}
(див. відповідні означення в \cite[п. 1.2]{Agranovich97}). Тут $L$~--- лінійний диференціальний оператор на $M$ довільного парного порядку $2r\geq2$, а кожне $B_{j}$~--- крайовий лінійний диференціальний оператор на $\partial M$ довільного порядку $m_{j}\geq0$. Усі коефіцієнти цих операторів належать до класів $C^{\infty}(M)$ і $C^{\infty}(\partial M)$ відповідно. Покладемо $B:=(B_{1},\ldots,B_{r})$ і $m:=\mathrm{max}\{m_{1},\ldots,m_{r}\}$.

Для кожного $j\in\{1,\ldots,r\}$ відображення $u\mapsto B_{j}u$, де $u\in C^{\infty}(M)$, продовжується єдиним чином (за неперервністю) до обмеженого лінійного оператора
\begin{equation}\label{bound-oper-H}
\begin{gathered}
B_{j}:B^{s}_{p,p}(M^{\circ})\to B^{s-m_j-1/p}_{p,p}(\partial M)\hookrightarrow\mathcal{D}'(\partial M)\\
\mbox{для усіх}\;\;p\in[1,\infty),\;s>m_j+1/p.
\end{gathered}
\end{equation}
Це випливає з \cite[с.~282, теорема п.~3.3.3]{Triebel86}. Позначимо через $\Upsilon_{m}(M^{\circ})$ об'єднання усіх просторів $B^{s}_{p,p}(M^{\circ})$, де $p\in[1,\infty)$ і $s>m+1/p$. З огляду на \eqref{bound-oper-H} і неперервність оператора $L$ на просторі $\overline{\mathcal{D}}\,'(M^{\circ})$ маємо коректно означений лінійний оператор
\begin{equation}\label{LB-general}
(L,B):\Upsilon_{m}(M^{\circ})\to \mathcal{D}\,'(M^{\circ})\times (\mathcal{D}'(\partial M))^{r}.
\end{equation}
Останній оператор буде неперервним, якщо на лінійному просторі $\Upsilon_{m}(M^{\circ})$ увести топологію індуктивної границі. У цей простір вкладені усі класичні простори $B^{s}_{p,q}(M^{\circ})$, де $p,q\in[1,\infty]$, і $F^{s}_{p,q}(M^{\circ})$, де $p,q\in[1,\infty)$, за умови $s>m+1/p$ для кожного з них. Це випливає з \cite[с.~278, теорема п.~3.3.1]{Triebel86}. Отже, у простір $\Upsilon_{m}(M^{\circ})$ вкладені узагальнені простори $B^{\alpha}_{p,q}(M^{\circ})$, де $p,q\in[1,\infty]$, і $F^{\alpha}_{p,q}(M^{\circ})$, де $p,q\in[1,\infty)$, за умови, що $\alpha\in\mathrm{RO}$ і $\sigma_{0}(\alpha)>m+1/p$. Тут, нагадаємо, $\sigma_{0}(\alpha)$~---
нижній індекс Матушевської функції~$\alpha$.

\smallskip

\textbf{Теорема 12.} \it Нехай $p,q\in[1,\infty]$, $\alpha\in\mathrm{RO}$ і $\sigma_{0}(\alpha)>m+1/p$. Тоді звуження відображення \eqref{LB-general} на простір $B^{\alpha}_{p,q}(M^{\circ})$ є нетеровим обмеженим оператором
\begin{equation}\label{LB-B}
(L,B):B^{\alpha}_{p,q}(M^{\circ})\to B^{\alpha\varrho^{-2r}}_{p,q}(M^{\circ})\times
\prod_{j=1}^{r}B^{\alpha\varrho^{-m_j-1/p}}_{p,q}(\partial M).
\end{equation}
Крім того, якщо $p,q<\infty$, то звуження відображення \eqref{LB-general} на простір $F^{\alpha}_{p,q}(M^{\circ})$ є нетеровим обмеженим оператором
\begin{equation}\label{LB-F}
(L,B):F^{\alpha}_{p,q}(M^{\circ})\to F^{\alpha\varrho^{-2r}}_{p,q}(M^{\circ})\times
\prod_{j=1}^{r}B^{\alpha\varrho^{-m_j-1/p}}_{p,p}(\partial M).
\end{equation}
Ядро кожного з цих операторів лежить в $C^{\infty}(M)$ і разом з індексом не залежить від параметрів $p$, $q$ і~$\alpha$. \rm

\smallskip

Вказані ядро та індекс називають також ядром та індексом досліджуваної ЕКЗ.

У випадку степеневої функції $\alpha$ ця теорема доведена в \cite[с.~82, теорема~5.2]{Johnsen96} (у більш загальній ситуації псевдодиференціальних ЕКЗ, які утворюють алгебру Буте де Монвеля). (Для регулярних диференціальних ЕКЗ, де $M^{\circ}$~--- евклідова область, ця теорема у вказаному випадку доведена трохи раніше в \cite[с.~145, теорема~14, с.~146, теорема~15]{FrankeRunst95}). У загальній ситуації теорема~12 виводиться з випадку степеневої функції $\alpha$ за допомогою інтерполяції на підставі теореми~9 і теореми про інтерполяцію нетерових операторів (остання аналогічна \cite[с.~35, теорема~1.7]{MikhailetsMurach10}). При цьому використовуються аналоги формул \eqref{B-gen-R-interp} і \eqref{F-gen-C-interp-F} для просторів на  $M^{\circ}$ і $\partial M$.

Версія цієї теореми правильна і для матричних ЕКЗ вигляду
\begin{equation}\label{matrix-bvp}
\mathbf{Lu}=\mathbf{f}\quad\mbox{на}\;\,M^\circ,\quad
\mathbf{Bu}=\mathbf{g}\quad\mbox{на}\;\,\partial M.
\end{equation}
Тут $\mathbf{L}=(L_{j,k})_{j,k=1}^{\lambda}$~--- матричний диференціальний оператор, еліптичний за Дуглісом--Ніренбергом на $M$ для набору цілих чисел $r_{1},\ldots,r_{\lambda}$ і $l_{1},\ldots,l_{\lambda}$; отже, кожне $\mathrm{ord}\,A_{j,k}\leq r_{j}+l_{k}$. Окрім того, $2r:=r_{1}+l_{1}+\cdots+r_{\lambda}+l_{\lambda}$ є парним числом, а $\mathbf{B}=(B_{j,k})_{j,k=1}^{r,\lambda}$~--- $r\times\lambda$-матриця скалярних крайових диференціальних операторів $B_{j,k}$, порядки яких задовольняють умову $\mathrm{ord}\,B_{j,k}\leq m_{j}+l_{k}$, де задано цілі числа $m_{1},\ldots,m_{r}$. Усі коефіцієнти скалярних диференціальних операторів $L_{j,k}$ і $B_{j,k}$ є нескінченно гладкими функціями на $M$ і $\partial M$ відповідно. Припускається, що ця матрична крайова задача є еліптичною за Агмоном--Дуглісом--Ніренбергом (див., наприклад, \cite[п.~6.1~a]{Agranovich97}). Для неї (як і у скалярному випадку) покладемо $m:=\max\{m_{1},\ldots,m_{r}\}$. З обмеженості операторів \eqref{bound-oper-H} випливає, що відображення
$\mathbf{u}\mapsto(\mathbf{L},\mathbf{B})\mathbf{u}$ означене коректно на парі просторів
\begin{equation}\label{LB-matrix-general}
(\mathbf{L},\mathbf{B}):\prod_{k=1}^{\lambda}\Upsilon_{l_k+m}(M^{\circ})\to
(\mathcal{D}\,'(M^{\circ}))^{\lambda}\times (\mathcal{D}'(\partial M))^{r}.
\end{equation}

\smallskip

\textbf{Теорема 12М.} \it Нехай $p,q\in[1,\infty]$, $\alpha\in\mathrm{RO}$ і $\sigma_{0}(\alpha)>m+1/p$. Тоді відповідне звуження відображення \eqref{LB-matrix-general} є нетеровим обмеженим оператором на парі просторів
\begin{equation}\label{LB-matrix-B}
(\mathbf{L},\mathbf{B}):
\prod_{k=1}^{\lambda}B^{\alpha\varrho^{l_k}}_{p,q}(M^{\circ})\to
    \prod_{j=1}^{\lambda}B^{\alpha\varrho^{-r_j}}_{p,q}(M^{\circ})\times
    \prod_{j=1}^{r}B^{\alpha\varrho^{-m_j-1/p}}_{p,q}(\partial M).
\end{equation}
Крім того, якщо $p,q<\infty$, то відповідне звуження відображення \eqref{LB-matrix-general} є нетеровим обмеженим оператором на парі просторів
\begin{equation}\label{LB-matrix-F}
(\mathbf{L},\mathbf{B}):
\prod_{k=1}^{\lambda}F^{\alpha\varrho^{l_k}}_{p,q}(M^{\circ})\to
    \prod_{j=1}^{\lambda}F^{\alpha\varrho^{-r_j}}_{p,q}(M^{\circ})\times
    \prod_{j=1}^{r}B^{\alpha\varrho^{-m_j-1/p}}_{p,p}(\partial M).
\end{equation}
Ядро кожного з цих операторів лежить в $(C^{\infty}(M))^{\lambda}$ і разом з індексом не залежить від параметрів $p$, $q$ і~$\alpha$. \rm

\smallskip

З теореми 12 і 12М випливають точні достатні умови локальної регулярності розв'язків еліптичних крайових задач в узагальнених $B$- і $F$-просторах та відповідні апріорні оцінки розв'язків (порівняти, наприклад, з \cite[сс. 301--302, теореми~3 і~4]{AnopKasirenko16MFAT}). Версії цих теорем правильні і для матричних псевдодиференціальних ЕКЗ, які утворюють алгебру Буте де Монвеля.

З крайовими операторами тісно пов'язане питання про точні простори слідів на $\partial M$ розподілів з узагальнених $B$- або $F$-просторів, заданих на $M^{\circ}$. Для узагальнених $B$-просторів і слідів на гіперплощинах воно було поставлено Ґ.~Трібелем в \cite[с.~57]{Triebel10}. Дамо відповідь на це питання. Нехай $\partial_{\nu}$ позначає оператор трансверсальної похідної функції, заданої на многовиді $M$. (Якщо $M$~--- евклідова область, то $\partial_{\nu}$~--- похідна уздовж внутрішньої нормалі до межі цієї області). Трансверсальна похідна означена у деякому околі межі $\partial M$. Для кожного $r\in\mathbb{N}$ розглянемо лінійне відображення
\begin{equation}\label{trace-r}
R_{r}:u\mapsto\bigl(u\!\upharpoonright\!\partial M,\ldots,
(\partial_{\nu}^{r-1}u)\!\upharpoonright\!\partial M\bigr),\quad
\mbox{де}\quad u\in\Upsilon_{r-1}(M^{\circ}).
\end{equation}

\textbf{Теорема 13.} \it Нехай $r\in\mathbb{N}$, $p,q\in[1,\infty]$, $\alpha\in\mathrm{RO}$ і $\sigma_{0}(\alpha)>r-1+1/p$. Тоді звуження відображення \eqref{trace-r} на простір $B^{\alpha}_{p,q}(M^{\circ})$ є  сюр'єктивним обмеженим оператором
\begin{equation}\label{trace-r-B}
R_{r}:B^{\alpha}_{p,q}(M^{\circ})\to
\prod_{j=1}^{r}B^{\alpha\varrho^{-j+1-1/p}}_{p,q}(\partial M).
\end{equation}
Крім того, якщо $p,q<\infty$, то звуження відображення \eqref{trace-r} на простір $F^{\alpha}_{p,q}(M^{\circ})$ є сюр'єктивним обмеженим оператором
\begin{equation}\label{trace-r-F}
R_{r}:F^{\alpha}_{p,q}(M^{\circ})\to
\prod_{j=1}^{r}B^{\alpha\varrho^{-j+1-1/p}}_{p,p}(\partial M).
\end{equation}
Більше того, існує лінійне відображення
\begin{equation}\label{T-r}
T_{r}:(\mathcal{D}'(\partial M))^{r}\to \overline{\mathcal{D}}\,'(M^{\circ})\cap C^{\infty}(M^{\circ}),
\end{equation}
незалежне від $p$, $q$ і $\alpha$, яке задає обмежені праві обернені оператори до \eqref{trace-r-B} і \eqref{trace-r-F}. А саме, відповідні звуження відображення \eqref{T-r} обмежені на парах просторів
\begin{gather}\label{T-r-B}
T_{r}:\prod_{j=1}^{r}B^{\alpha\varrho^{-j+1-1/p}}_{p,q}(\partial M)\to
B^{\alpha}_{p,q}(M^{\circ}),\\
T_{r}:\prod_{j=1}^{r}B^{\alpha\varrho^{-j+1-1/p}}_{p,p}(\partial M)\to
F^{\alpha}_{p,q}(M^{\circ})\label{T-r-F}
\end{gather}
та задовольняють умову
\begin{equation}\label{RrTr}
R_{r}T_{r}g=g\quad\mbox{для довільного}\quad
g\in\prod_{j=1}^{r}B^{\alpha\varrho^{-j+1-1/p}}_{p,q}(\partial M).
\end{equation}
\rm

Обговоримо доведення цієї теореми. Нехай $L$~--- формально самоспряжений (відносно скалярного добутку в $L_{2}(M,dx)$) правильно еліптичний диференціальний оператор на $M$ порядку $2r$ і з коефіцієнтами класу $C^{\infty}(M)$. Розглянемо регулярну ЕКЗ
\begin{gather}
(L-i)u=f\quad\mbox{на}\;\,M^\circ,\\
\partial_{\nu}^{j-1}u=g_{j}\quad\mbox{на}\;\,\partial M,
\quad j=1,...,r.
\end{gather}
Згідно з \cite[c.~477, теорема п.~5.4.2]{Triebel80} і \cite[c.~491, теорема (b) п.~5.5.2]{Triebel80} ця задача має нульові ядро та індекс. Тому за теоремою~12 маємо ізоморфізми
\begin{gather}\label{L-iI-R-B}
(L-iI,R_{r}):B^{\alpha}_{p,q}(M^{\circ})\leftrightarrow B^{\alpha\varrho^{-2r}}_{p,q}(M^{\circ})\times
\prod_{j=1}^{r}B^{\alpha\varrho^{-j+1-1/p}}_{p,q}(\partial M),\\
(L-iI,R_{r}):F^{\alpha}_{p,q}(M^{\circ})\leftrightarrow F^{\alpha\varrho^{-2r}}_{p,q}(M^{\circ})\times
\prod_{j=1}^{r}B^{\alpha\varrho^{-j+1-1/p}}_{p,p}(\partial M),
\end{gather}
де $I$~--- тотожний оператор на $\overline{\mathcal{D}}\,'(M^{\circ})$. Звідси негайно випливає обмеженість і сюр'єктивність операторів \eqref{trace-r-B} і \eqref{trace-r-F}. Уведемо відображення \eqref{T-r}. Відображення $u\mapsto R_{r}u$, задане на усіх функціях $u\in C^{\infty}(M)$ таких, що $(L-i)u=0$ на $M$, продовжується за неперервністю до ізоморфізмів
\begin{equation}
R_{r}:H^{s}_{2}(M^{\circ},L-iI)\leftrightarrow
\prod_{j=1}^{r}H^{s-j+1/2}_{2}(\partial M)
\quad\mbox{для усіх}\quad s\in\mathbb{R}
\end{equation}
(див., наприклад, \cite[с.~165, теорема~3.11]{MikhailetsMurach10} у випадку, коли $M^{\circ}$~--- евклідова область). Тут
\begin{equation}
H^{s}_{2}(M^{\circ},L-iI):=\{u\in H^{s}_{2}(M^{\circ}):(L-i)u=0\;\mbox{на}\;M^{\circ}\}
\end{equation}
--- підпростір гільбертового простору $H^{s}_{2}(M^{\circ})$. Обернені оператори до цих ізоморфізмів коректно визначають лінійне відображення \eqref{T-r} з огляду на рівність
\begin{equation}
\mathcal{D}'(\partial M)=\bigcup_{l\in\mathbb{R}}H^{l}_{2}(\partial M)
\end{equation}
і включення
\begin{equation}
H^{s}_{2}(M^{\circ},L-iI)\subset C^{\infty}(M^{\circ}).
\end{equation}
Звуження відображення \eqref{T-r} на простір
\begin{equation}
\prod_{j=1}^{r}B^{\alpha\varrho^{-j+1-1/p}}_{p,q}(\partial M)
\end{equation}
збігається з оператором, оберненим до звуження ізоморфізму \eqref{L-iI-R-B} на підпростір
\begin{equation}
\{u\in B^{\alpha}_{p,q}(M^{\circ}):
(L-i)u=0\;\mbox{на}\;M^{\circ}\}
\end{equation}
простору $B^{\alpha}_{p,q}(M^{\circ})$. Це тягне за собою обмеженість оператора \eqref{T-r-B} і рівність \eqref{RrTr}. Аналогічно міркуємо і для оператора \eqref{T-r-F}.

З теореми~13 при $r=1$ випливає такий опис “позитивних” узагальнених $B$-просторів на $\partial M$ у термінах слідів (порівняти з \cite[сс. 155--156, доведення наслідку 3.1]{MikhailetsMurach10}).

\smallskip

\textbf{Теорема 14.} \it Нехай $p,q\in[1,\infty]$, $\beta\in\mathrm{RO}$ і $\sigma_{0}(\beta)>0$. Тоді
\begin{gather}
B^{\beta}_{p,q}(\partial M)=\bigl\{h:=u\!\upharpoonright\!\partial M:
u\in B^{\beta\varrho^{1/p}}_{p,q}(M^{\circ})\bigr\},\\
\|h\|_{B^{\beta}_{p,q}(\partial M)}\asymp
\inf\bigl\{\|u\|_{B^{\beta\varrho^{1/p}}_{p,q}(M^{\circ})}:
u\!\upharpoonright\!\partial M=h\bigr\}.
\end{gather}
Крім того, якщо $p,q\neq\infty$, то
\begin{gather}
B^{\beta}_{p,p}(\partial M)=\bigl\{h:=u\!\upharpoonright\!\partial M:
u\in F^{\beta\varrho^{1/p}}_{p,q}(M^{\circ})\bigr\},\\
\|h\|_{B^{\beta}_{p,p}(\partial M)}\asymp
\inf\bigl\{\|u\|_{F^{\beta\varrho^{1/p}}_{p,q}(M^{\circ})}:
u\!\upharpoonright\!\partial M=h\bigr\}.
\end{gather} \rm

Тут, як звичайно, символ $\asymp$ позначає еквівалентність норм.

\medskip

\textbf{7. Квазібанахові простори узагальненої гладкості.} Означення 3 і 4 просторів $B^{\alpha}_{p,q}(\mathbb{R}^{n})$ і $F^{\alpha}_{p,q}(\mathbb{R}^{n})$, де $\alpha\in\mathrm{RO}$, є змістовними і у випадку, коли $p,q\in(0,\infty]$ для $B$-просторів, або коли $p,q\in(0,\infty)$ для $F$-просторів. Якщо $0<p<1$ та/або $0<q<1$, то отримаємо квазібанахові простори. (Означення останніх відрізняються від означення банахових просторів лише тим, що замість нерівності трикутника для норм виконується більш слабка нерівність $\|u+v\|\leq c\,(\|u\|+\|v\|)$ з деякою сталою $c\geq1$.) У випадку, коли $\alpha(t)\equiv t^{s}$ для деякого $s\in\mathbb{R}$, отримаємо класичні квазібанахові простори Бєсова $B^{s}_{p,q}(\mathbb{R}^{n})$ і Лізоркіна--Трібеля $F^{s}_{p,q}(\mathbb{R}^{n})$. Їх досліджували і застосовували у багатьох роботах; див. монографії Ґ.~Трібеля \cite{Triebel86, Triebel06, Triebel92} та наведені там посилання.

Метод дійсної інтерполяції з функціональним параметром $\gamma\in\mathfrak{B}$, підпорядкованим умові
\eqref{Boyd-indexes}, та числовим параметром $q\in(0,\infty]$ є інтерполяційним функтором і у випадку квазібанахових просторів \cite[с.~205]{Persson86}.

\smallskip

\textbf{Теорема 15.} \it Твердження теорем $4$ і $5$ є правильними, якщо $q\in(0,\infty]$, а $p,q_{0},q_{1}\in(0,\infty]$ для $B$-просторів та
$p,q_{0},q_{1}\in(0,\infty)$ для $F$-просторів. \rm

\smallskip

Теорема 15 доводиться у спосіб, використаний у статті \cite[c.~166]{CobosFernandez88} для банахових просторів узагальненої гладкості; це зазначено у праці \cite[п.~4.2, с.~705]{AlmeidaCaetano11}. Частина цієї теореми, яка стосується формули \eqref{B-gen-R-interp}, наведена у статтях \cite[с.~205, твердження~7]{Almeida05} і \cite[п.~5.1, с.~708]{AlmeidaCaetano11}. Згідно з цією формулою кожний простір $B^{\alpha}_{p,q}(\mathbb{R}^{n})$, де $\alpha\in\mathrm{RO}$ і $p,q\in(0,\infty]$, можна отримати за допомогою дійсної інтерполяції класичних просторів Нікольського\,--\,Бєсова.

Аналог останньої властивості для $F$-просторів є правильним щодо інтерполяційного методу Густавсона\,--\,Петре, відомого також як $\pm$-метод інтерполяції (див. його означення в \cite[п.~6]{GustavssonPeetre77} та \cite[п.~1]{Gustavsson82}). Для довільної інтерполяційної пари квазібанахових просторів $E_0$ і $E_1$ позначимо через $\langle E_0,E_1\rangle_{\gamma}$ квазібанахів простір, отриманий інтерполяцією цієї пари вказаним методом, тут параметром інтерполяції є довільна функція $\gamma\in\mathfrak{B}$, яка задовольняє умову \eqref{Boyd-indexes}. Цей метод є інтерполяційним функтором, який має властивість, аналогічну до \eqref{R-interp-property}, тобто
\begin{equation}
\|L\|_{\langle E_0,E_1\rangle_{\gamma}\to\langle Q_{0},Q_{1}\rangle_{\gamma}}\leq
c\cdot\|L\|_{E_{0}\to Q_{0}}\cdot\overline{\gamma}
\biggl(\frac{\|L\|_{E_{1}\to Q_{1}}}{\|L\|_{E_{0}\to Q_{0}}}\biggr),
\end{equation}
де число $c>0$ не залежить від оператора $L$ та пар просторів $E_{0},E_{1}$ і $Q_{0},Q_{1}$ \cite[твердження 6.1]{GustavssonPeetre77}.

\smallskip

\textbf{Теорема 16.} \it Нехай задано показники регулярності  $\alpha_{0},\alpha_{1}\in\mathrm{RO}$, показники сумовності  $p,q\in(0,\infty)$ та параметр інтерполяції $\gamma\in\mathfrak{B}$, який задовольняє умову \eqref{Boyd-indexes}. Означимо параметр $\alpha\in\mathrm{RO}$ за формулою \eqref{alpha-R-interp-gamma}. Тоді
\begin{equation}\label{GP-interp-F-gen}
\langle F^{\alpha_0}_{p,q}(\mathbb{R}^{n}),
F^{\alpha_1}_{p,q}(\mathbb{R}^{n})\rangle_{\gamma}=
F^{\alpha}_{p,q}(\mathbb{R}^{n})
\end{equation}
з точністю до еквівалентності норм. \rm

\smallskip

Ця теорема доводиться за допомогою методики, розробленої в \cite[п.~12]{FrazierJawerth90} для класичних просторів Лізоркіна\,--\,Трібеля.

З теореми~16 випливає, що кожний простір  $F^{\alpha}_{p,q}(\mathbb{R}^{n})$, де $\alpha\in\mathrm{RO}$, можна отримати, інтерполюючи методом Густавсона\,--\,Петре деякі пари класичних просторів Лізоркіна\,--\,Трібеля.

\smallskip

\textbf{Теорема 17.} \it Нехай $\alpha\in\mathrm{RO}$ і $p,q\in(0,\infty)$. Довільно виберемо дійсні числа $\nobreak{s_0<\sigma_0(\alpha)}$ і $s_1>\sigma_1(\alpha)$ та означимо інтерполяційний параметр $\gamma$ так само, як у теоремі~$5$. Тоді
\begin{equation}\label{F-gen-GP-interp}
F^{\alpha}_{p,q}(\mathbb{R}^{n})=
\langle F^{s_0}_{p,q}(\mathbb{R}^{n}),
F^{s_1}_{p,q}(\mathbb{R}^{n})\rangle_{\gamma}
\end{equation}
з точністю до еквівалентності норм. \rm

За допомогою міркувань, наведених у п.~4, доводяться аналоги цих теорем для квазінормованих узагальнених просторів Бєсова і Лізоркіна\,--\,Трібеля, заданих на деяких евклідових областях або на довільному компактному орієнтовному $C^{\infty}$-многовиді $M$ (з краєм чи без).

\smallskip

\textbf{Теорема 18.} Твердження теорем $15$--$17$ є правильними, якщо у них замінити $\mathbb{R}^{n}$ на півпростір $\mathbb{R}^{n}_{+}$ або на довільну обмежену область $G\subset\mathbb{R}^{n}$ з ліпшіцевою межею, або на $M^{\circ}:=M\setminus\partial M$. \rm

\smallskip

За допомогою теореми 18 доводяться версії усіх інших (неінтерполяційних) результатів пп. 4--6 для вказаних просторів. Сформулюємо відповідні теореми. Розглянемо умову
\begin{equation}\label{assumption-p-q}
\left\{
  \begin{array}{ll}
    p,q\in(0,\infty]&\hbox{у випадку $B$-просторів;} \\
    p,q\in(0,\infty)&\hbox{у випадку $F$-просторів.}
  \end{array}
\right.
\end{equation}

\smallskip

\textbf{Теорема 19.} \it Властивості універсального лінійного оператора продовження, сформульовані у наслідку~$3$, виконуються для довільного $\alpha\in\mathrm{RO}$ за умови \eqref{assumption-p-q}. \rm

\smallskip

\textbf{Теорема 20.} \it Твердження теорем $8$ і $10$ є правильними для довільного $\alpha\in\mathrm{RO}$ за умови \eqref{assumption-p-q}. \rm

\smallskip

Наступні дві теореми стосуються застосувань квазібанахових узагальнених просторів Бєсова і Лізоркіна\,--\,Трібеля на многовидах до еліптичних операторів і еліптичних крайових задач.

\smallskip

\textbf{Теорема 21.} \it Нехай $\partial M=\varnothing$, та $A$ і $\mathbf{A}$ є еліптичні ПДО, що фігурують у теоремах \rm 11 \it і \rm 11М\it, відповідно. Тоді правильні такі твердження:
\begin{itemize}
  \item[\rm (i)\rm] Якщо простори  \eqref{kernel-A} і \eqref{cokernel-A} нульові, то оператори \eqref{PsDO-B} і \eqref{PsDO-F} є топологічними ізоморфізмами для довільного $\alpha\in\mathrm{RO}$ за умови \eqref{assumption-p-q}.
  \item[\rm (ii)\rm] Якщо ядро і коядро ПДО $\mathbf{A}$ нульові, то
  оператори \eqref{PsDO-matrix-B} і \eqref{PsDO-matrix-F} є топологічними ізоморфізмами для довільного $\alpha\in\mathrm{RO}$ за умови \eqref{assumption-p-q}.
\end{itemize} \rm

\smallskip

Розглянемо еліптичні крайові задачі, скалярну \eqref{ep-equ}, \eqref{ep-bound-cond} і матричну \eqref{matrix-bvp}, де $n:=\dim M\geq2$ і $\partial M\neq\varnothing$. У цьому підрозділі позначимо через $\Upsilon_{\ell}(M^{\circ})$, де $\ell\in\mathbb{R}$, об'єднання усіх просторів $B^{s}_{p,p}(M^{\circ})$ таких, що $p\in(0,\infty)$ і $s>\ell+\lambda(p,n)$, де
\begin{equation*}
\lambda(p,n):=\frac{1}{p}+
\max\biggl\{0,(n-1)\biggl(\frac{1}{p}-1\biggl)\biggr\}.
\end{equation*}
Згідно з \cite[с.~193, теорема п.~2.7.2]{Triebel86} означені коректно  лінійні оператори \eqref{LB-general} і \eqref{LB-matrix-general}, де для обох задач $m:=\max\{m_{1},\ldots,m_{r}\}$. Звісно, $\lambda(p,n)=1/p$, якщо $1\leq p\leq\infty$. 

\smallskip

\textbf{Теорема 22.} \it Припустимо, що параметри $\alpha\in\mathrm{RO}$, $p$ та $q$, задовольняють умови \eqref{assumption-p-q} і $\nobreak{\sigma_{0}(\alpha)>m+\lambda(p,n)}$. Тоді правильні такі твердження:
\begin{itemize}
  \item[\rm (i)\rm] Якщо ядро та індекс ЕКЗ \eqref{ep-equ}, \eqref{ep-bound-cond} нульові, то відповідні звуження відображення \eqref{LB-general} встановлюють топологічні ізоморфізми на парах просторів, вказаних у формулах \eqref{LB-B} і \eqref{LB-F}.
  \item[\rm (ii)\rm] Якщо ядро та індекс ЕКЗ \eqref{matrix-bvp} нульові, то відповідні звуження відображення \eqref{LB-matrix-general} встановлюють топологічні ізоморфізми на парах просторів, вказаних у формулах \eqref{LB-matrix-B} і \eqref{LB-matrix-F}.
\end{itemize} \rm

\smallskip

У загальному випадку ненульових ядер та індексів еліптичних ПДО і еліптичних крайових задач мають місце версії теорем 21 і 22, сформульовані у термінах, використаних у \cite[с.~82, теорема~5.2]{Johnsen96} і \cite[с.~145, теорема~14; с.~146, теорема~15]{FrankeRunst95}.

На завершення сформулюємо версії теорем 13 і 14 про точні простори слідів на межі $\partial M$ для розглянутих у цьому підрозділі квазібанахових просторів узагальненої гладкості.

\smallskip

\textbf{Теорема 23.} \it Нехай $r\in\mathbb{N}$. Припустимо, що параметри $\alpha\in\mathrm{RO}$, $p$ і $q$ задовольняють умови \eqref{assumption-p-q} і $\sigma_{0}(\alpha)>r-1+\lambda(p,n)$. Тоді є правильним висновок теореми~$13$. \rm

\smallskip

\textbf{Теорема 24.} \it Припустимо, що параметри $\beta\in\mathrm{RO}$, $p$ і $q$ задовольняють умови \eqref{assumption-p-q} і $\sigma_{0}(\beta)>\lambda(p,n)-1/p$. Тоді є правильним висновок теореми~$14$. \rm

\end{document}